\input ./style/arxiv-general.cfg
\documentclass[MSNbibl,number,citesort,seceqn,dvips]{arxbj}
\makeatletter
   \@ifpackageloaded{graphicx}{}{\usepackage{graphicx}}
\makeatother
\usepackage{upgreek}

\volume{22}
\issue{3}
\pubyear{2016}
\firstpage{1770}
\lastpage{1807}
\doi{10.3150/15-BEJ711}
\docsubty{FLA}

\makeatletter
\newcommand{\upi}{\uppi}
\newcommand{\dd}{\mathrm{d}}
\newcommand{\rrvert}{\vert}
\newcommand{\rrVert}{\Vert}
\newcommand{\llvert}{\vert}
\newcommand{\llVert}{\Vert}
\renewcommand{\mid}{|}
\newcommand{\xrightarrow}[1]{\stackrel{#1}{\to}}
\newcommand{\mod}{\quad\mathrm{mod}~}
\newcommand{\modt}{~~\mathrm{mod}~}
\newcommand{\mathds}{\mathbb}
\newcommand{\IE}{\mathds{E}} 
\newcommand{\IG}{\mathds{G}}
\newcommand{\IN}{\mathds{N}}
\newcommand{\IP}{\mathds{P}}
\newcommand{\IR}{\mathds{R}}
\newcommand{\IZ}{\mathds{Z}}
\newcommand{\ee}{\mathrm{e}}
\newcommand{\ii}{\mathrm{i}}
\newcommand{\eps}{\varepsilon}
\newcommand{\E}{\mathds{E}}
\newcommand{\cF}{\mathcal{F}}
\newcommand{\Cov}{\operatorname{Cov}}
\newcommand{\Var}{\operatorname{Var}}
\newcommand{\cum}{\operatorname{cum}}
\newcommand{\pkg}[1]{\textbf{#1}}
\newcommand{\proglang}{\textsf}
\newtheorem{theorem}{Theorem}[section]
\newtheorem{lemma}[theorem]{Lemma}
\newremark{bem}{Remark}[section]
\newremark{re}{Remark}
\newtheorem{prop}[theorem]{Proposition}
\newcommand{\PerProc}{\mathbb{I}}
\newcommand{\weak}{\rightsquigarrow}
\newcommand{\inr}{I_{n,R}}
\makeatother

\begin{document}
\begin{frontmatter}

\title{Quantile spectral processes: Asymptotic~analysis and inference}
\runtitle{Quantile spectral processes}

\begin{aug}
\author[A]{\inits{T.}\fnms{Tobias}~\snm{Kley}\thanksref{A,e3}\ead[label=e3,mark]{tobias.kley@rub.de}},
\author[A]{\inits{S.}\fnms{Stanislav}~\snm{Volgushev}\thanksref{A,e4}\ead[label=e4,mark]{stanislav.volgushev@rub.de}},
\author[A]{\inits{H.}\fnms{Holger}~\snm{Dette}\corref{}\thanksref{A,e1}\ead[label=e1,mark]{holger.dette@rub.de}}
\and
\author[B,C]{\inits{M.}\fnms{Marc}~\snm{Hallin}\thanksref{B,C,e2}\ead[label=e2,mark]{mhallin@ulb.ac.be}}
\address[A]{Fakult\"{a}t f\"{u}r Mathematik, Ruhr-Universit\"{a}t
Bochum, 44780 Bochum, Germany.\\ \printead{e3,e4,e1}}
\address[B]{ECARES, Universit\'{e} Libre de Bruxelles, 1050 Brussels,
Belgium. \printead{e2}}
\address[C]{ORFE, Princeton University, Princeton, NJ 08544, USA}
\end{aug}

%
\received{\smonth{7} \syear{2014}}
%
\revised{\smonth{2} \syear{2015}}

%
\begin{abstract}
Quantile- and copula-related spectral concepts recently have been
considered by various authors. Those spectra, in their most general
form, provide a full characterization of the copulas associated with
the pairs $(X_t, X_{t-k})$ in a process $(X_t)_{t\in\mathbb{Z}}$, and
account for important dynamic features, such as changes in the
conditional shape (skewness, kurtosis), time-irreversibility, or
dependence in the extremes that their traditional counterparts cannot
capture. Despite various proposals for estimation strategies, only
quite incomplete asymptotic distributional results are available so far
for the proposed estimators, which constitutes an important obstacle
for their practical application. In this paper, we provide a detailed
asymptotic analysis of a class of smoothed rank-based
cross-periodograms associated with the \textit{copula spectral density
kernels} introduced in Dette \textit{et~al.}
[\textit{Bernoulli} \textbf{21} (2015) 781--831]. We show that, for a very
general class of (possibly nonlinear) processes, properly scaled and
centered smoothed versions of those cross-periodograms, indexed by
couples of quantile levels, converge weakly, as stochastic processes,
to Gaussian {processes}. A first application of those results is the
construction of asymptotic confidence intervals for \textit{copula
spectral density kernels}. The same convergence results also provide
asymptotic distributions (under serially dependent observations) for a
new class of rank-based spectral methods involving the Fourier
transforms of rank-based serial statistics such as the \textit
{Spearman, Blomqvist} or \textit{Gini autocovariance coefficients}.
\end{abstract}

%
\end{frontmatter}

\section{Introduction}\label{sec1}

Spectral analysis and frequency domain methods play a central role in
the nonparametric analysis of time series data. The classical frequency
domain representation is based on the \textit{spectral density}~--
call it the $L^2$-\textit{spectral density} in order to distinguish it
from other \textit{spectral densities} to be defined in the sequel --
which is traditionally defined as the Fourier transform of the
autocovariance function of the process under study. Fundamental tools
for the estimation of spectral densities are the \textit{periodogram} and
its \textit{smoothed} versions. The classical periodogram~-- similarly
call it the $L^2$-\textit{periodogram}~-- can be defined either as the
discrete Fourier transform of empirical autocovariances, or through
$L^2$-projections of the observed series on a harmonic basis. The
success of periodograms in time series analysis is rooted in their fast
and simple computation (through the \textit{fast Fourier transform}
algorithm) and their interpretation in terms of cyclic behavior, both
of a stochastic and of deterministic nature, which in specific
applications are more illuminating than time-domain representations.
$L^2$-periodograms are particularly attractive in the analysis of
Gaussian time series, since the distribution of a Gaussian process is
completely characterized by its spectral density. Classical references are
Priestley \cite{Priestley1981}, Brillinger \cite{Brillinger1975} or
Chapters~4 and~10
of Brockwell and Davis \cite{BrockwellDavis1987}.

Being intrinsically connected to means and covariances, the
$L^2$-spectral density and $L^2$-pe\-riodogram inherit the nice features
(such as optimality properties in the analysis of Gaussian series) of
$L^2$-methods, but also their weaknesses: they are lacking robustness
against outliers and heavy tails, and are
unable to capture important dynamic features such as changes in the
conditional shape (skewness, kurtosis), time-irreversibility, or
dependence in the extremes. This was realized by many researchers, and
various extensions and modifications of the $L^2$-periodogram have been
proposed to remedy those drawbacks.

Robust nonparametric approaches to frequency domain estimation have
been considered first; see
Kleiner, Martin and Thomson
\cite{KleinerEtAl1979} for an early
contribution, and Chapter~8 of Maronna, Martin and
Yohai \cite{MaronnaEtAl2006} for an overview.
More recently, Kl{\"u}ppelberg and Mikosch \cite
{KluppelbergMikosch1994} proposed a weighted
(``self-normalized'') version of the periodogram; see also Mikosch \cite
{Mikosch1998}. Hill and McCloskey \cite{HillMcCloskey2013} used a robust version of
autocovariances and a robustified periodogram with the goal to obtain
$L^2$-spectrum-based parameter estimates that are robust to
heavy-tailed data.
In the context of signal detection, Katkovnik
\cite{Katkovnik1998} introduced a
periodogram based on robust loss functions. The objective of all those
attempts is a robustification of classical tools: they essentially aim
at protecting existing $L^2$-spectral methods against the impact of
possible outliers or violations of distributional assumptions.

Other attempts, more recent and somewhat less developed, are
introducing alternative spectral concepts and tools, mostly related
with quantiles or copulas, and accounting for more general dynamic
features. A first step in that direction was taken by Hong \cite
{Hong1999}, who proposes a generalized spectral
density with
covariances replaced by joint characteristic functions. In the specific
problem of testing pairwise independence, Hong
\cite{Hong2000} introduces
a test statistic based on the Fourier transforms of (empirical) joint
distribution functions and copulas at different lags. Recently, there
has been a renewed surge of interest in that type of concept, with the
introduction, under the names of \textit{Lap\-lace-, quantile-} and
\textit{copula} spectral density and spectral density \textit{kernels},
of various quantile-related spectral concepts, along with the
corresponding sample-based periodograms and smoothed periodograms.
That strand of literature includes Li \cite{Li2008,Li2012,Li2013},
Hagemann \cite{Hagemann2013}, Dette \textit
{et~al.} \cite{DetteEtAl2013} and Lee and Rao
\cite{LeeRao2012}.
A Fourier analysis of extreme events, which is related in spirit but
quite different in many other respects, was considered by
Davis, Mikosch and Zhao
\cite{DavisMikoschZhao2013}. Finally, in
the time domain, Linton and Whang \cite
{LintonWhang2007}, Davis and Mikosch \cite
{DavisMikosch2009} and Han \textit{et~al.} \cite
{HanLintonOkaWhang2013} introduced the related concepts of \textit
{quantilograms} and \textit{extremograms}. A more detailed account of
some of those contributions is given in Section~\ref{sec2}.

A deep understanding of the distributional properties of any
statistical tool is crucial for its successful application. The
construction of confidence intervals, testing procedures and efficient
\mbox{estimators} all rest on results concerning finite-sample or asymptotic
properties of related statistics~-- here the appropriate (smoothed)
periodograms associated with the quantile-related spectral density
under study. Obtaining such asymptotic results, unfortunately, is not
trivial, and to the best of our knowledge, no results on the asymptotic
distribution of the aforementioned (smoothed) quantile and copula
periodograms are available so far.

In the case of i.i.d. observations, Hong \cite
{Hong2000} derived the
asymptotic distribution of an empirical version of the integrated
version of his quantile spectral density, while Lee and Rao \cite{LeeRao2012}
investigated the distributions of Cram\' er--von Mises-type statistics
based on empirical joint distributions. No results on the asymptotic
distribution of the periodogram itself are given, though. Li \cite{Li2008,Li2012} does not consider asymptotics
for smoothed versions of
his quantile periodograms, while the asymptotic results in Hagemann
\cite{Hagemann2013} and Dette \textit
{et~al.} \cite{DetteEtAl2013} are quite incomplete.
This is perhaps not so surprising: the asymptotic distribution of
classical $L^2$-spectral density estimators for general nonlinear
processes also has remained an active domain of research for several
decades; see Brillinger \cite{Brillinger1975}
for early results, Shao and Wu \cite{ShaoWu2007},
Liu and Wu \cite{LiuWu2010} or Giraitis and Koul \cite
{GiraitisKoul2013} for more
recent references.

The present paper has two major objectives. First, it aims at providing
a rigorous analysis of the asymptotic properties of a general class of
smoothed rank-based copula cross-periodograms generalizing the quantile
periodograms introduced by Hagemann \cite
{Hagemann2013} and, in an integrated
version, by Hong \cite{Hong2000}. In
Section~\ref{secasy}, we show that,
for general nonlinear time series, properly centered and smoothed
versions of those cross-periodograms, indexed by couples of quantile
levels, converge in distribution to centered Gaussian {processes}. A
first application of those results is the construction of asymptotic
confidence intervals which we discuss in detail in Section~\ref{secsim}.

The second objective of this paper is to introduce a new class of
rank-based frequency domain methods that can be described as a
non-standard rank-based Fourier analysis of the serial features of time
series. Examples of such methods are discussed in detail in
Section~\ref{secrankper}, where we study a class of
spectral densities, such as the \textit{Spearman}, \textit{Blomqvist and Gini
spectra}, and the corresponding periodograms, associated with
rank-based autocovariance concepts.
Denoting by $F$ the marginal distribution function of~$X_t$, the
\textit{Spearman spectral density}, for instance, is defined as $\sum_{k \in\IZ} \ee^{\ii\omega k}\rho_k^{\mathrm{Sp}}$,
where~$\rho_k^{\mathrm{Sp}}:=\operatorname{Corr}(F(X_t), F(X_{t-k}))$
denotes the lag-$k$ \textit{Spearman autocorrelation}. We show that
estimators of those spectral densities can be obtained as functionals
of the rank-based copula periodograms investigated in this paper. This
connection, and our process-level convergence results on the rank-based
copula periodograms, allow us to establish the asymptotic normality of
the smoothed versions of the newly defined rank-based periodograms.
Those results can be considered as frequency domain versions of H\'
ajek's celebrated asymptotic representation and normality results for
(non-serial) linear rank statistics under non-i.i.d. observations
(H{\'a}jek \cite{Hajek1968}).

The paper is organized as follows. Section~\ref{secdef} provides precise
definitions of the spectral concepts to be considered throughout, and
motivates the use of our quantile-related methods by a graphical comparison
of the copula spectra of white noise, QAR(1) and ARCH(1) processes,
respectively~-- all of which share the same helplessly flat $L^2$-spectral
density. Section~\ref{secasy} is devoted to the asymptotics of
rank-based copula (cross-)periodograms and their smoothed versions,
presenting the
main results of this paper: the convergence, for fixed frequencies
$\omega$, of the
smoothed copula rank-based periodogram indexed by couples of $(\tau_1,
\tau_2)$ of quantile orders to a Gaussian process
(Theorem~\ref{thmAsympDensityRankEstimator}). That theorem is based
on an equally
interesting asymptotic representation result
(Theorem~\ref{thmAsympDensityEstimator}). Section~\ref{secrankper}
deals with the relation with Spearman, Gini, and Blomqvist autocorrelation
coefficients and the related spectra. Based on a short Monte-Carlo study,
Section~\ref{secsim} discusses the practical performances of the methods
proposed, and Section~\ref{secConcl} provides some conclusions and
directions for future research. Proofs are concentrated in an \hyperref[pmt]{Appendix}
and an the online supplement~\cite{supp}.

\section{Copula spectral density kernels and rank-based periodograms}\label{secdef}\label{sec2}
In this section, we provide more precise
definitions of the various quantile- and copula-related spectra
mentioned in the \hyperref[sec1]{Introduction}, along with the corresponding periodograms.

Throughout, let $(X_t)_{t \in\IZ}$ denote a strictly stationary
process, of which we observe a finite stretch $X_0,\ldots,X_{n-1}$, say.
Denote by $F$ the marginal distribution function of $X_t$, and by
$q_\tau:= \inf\{x \in\IR\dvt  \tau\leq F(x)\}$, $\tau\in[0,1]$ the
corresponding quantile function, where we use the convention $\inf \varnothing =\infty$.
Note that if $\tau\in\{0,1\}$ then
$-\infty$ and $\infty$ are possible values for $q_{\tau}$. Our main
object of interest is the \textit{copula spectral density kernel}
\begin{equation}
\label{fqq1} \mathfrak{f}_{q_{\tau_1}, q_{\tau_2}}(\omega):= \frac{1}{2\uppi} \sum
_{k\in\IZ} \mathrm{e}^{-\ii\omega k}
\gamma^{U}_k(\tau _1,\tau _2),
\qquad\omega\in\IR, (\tau_1, \tau_2 )
\in[0,1]^2,
\end{equation}
based on the \textit{copula cross-covariances}
\[
\gamma^{U}_k(\tau_1,\tau_2):=
\textrm{Cov} \bigl(I\{ U_t\leq\tau_1\}, I\{
U_{t-k}\leq\tau_2\} \bigr), \qquad k\in\mathbb{Z},
\]
where $U_t:=F(X_t)$. Those copula spectral density kernels were
introduced in Dette \textit{et~al.} \cite{DetteEtAl2013},
and generalize the \textit{$\tau
$th quantile spectral densities} of Hagemann
\cite{Hagemann2013}, with which
they coincide for $\tau_1 = \tau_2 = \tau$; an integrated version
actually was first considered by Hong \cite{Hong2000}.
The same copula spectral density kernel also takes the
form
%
\begin{eqnarray}
\label{fqq2}
&& \mathfrak{f}_{q_{\tau_1}, q_{\tau_2}}(\omega)
\nonumber\\[-8pt]\\[-8pt]\nonumber
&&\quad := \frac{1}{2\uppi} \sum
_{k\in\IZ} \mathrm{e}^{-\ii\omega k} \bigl(
\IP(X_k \leq q_{\tau
_1},X_{ 0 } \leq
q_{\tau_2}) - \tau_1\tau_2 \bigr),
\qquad
\omega\in\IR, (\tau_1,\tau_2)
\in[0,1]^2,
\end{eqnarray}
%
where $\IP(X_k\leq q_{\tau_1}, X_0 \leq q_{\tau_2})$ is the joint distribution function of the pair $(X_k,X_0)$
taken at $(q_{\tau_1}, q_{\tau_2})$. This is, by definition, the copula of the pair $(X_k,X_0)$ evaluated
at $(\tau_1, \tau_2)$, while $\tau_1\tau_2$ is the independence copula evaluated at the same $(\tau_1, \tau_2)$.
The copula spectral density kernel thus can be interpreted as the
Fourier transform of the difference between pairwise copulas at lag $k$
and the independence copula, which justifies the notation and the terminology.

As shown by Dette \textit{et~al.} \cite{DetteEtAl2013}, the
copula spectral densities
provide a complete description of the pairwise copulas of a time
series. Similar to the regression setting, where joint distributions
and quantiles provide more information than covariances and means, the
copula spectral density kernel, by accounting for much more than the
covariance structure of a series, extends and supplements the classical
$L^2$-spectral density.

%
%
\begin{figure}[b]

\includegraphics{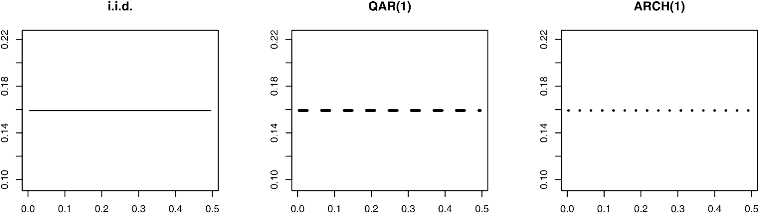}

\caption{Traditional $L^2$-spectra $(2\uppi)^{-1} \sum_{k \in\IZ}
\Cov(Y_{t+k}, Y_t) \mathrm{e}^{-\mathrm{i} \omega k}$. The process
$(Y_t)$ in the left-hand picture is independent standard normal white
noise; in the middle picture, $Y_t = X_t / \Var(X_t)^{1/2}$ where
$(X_t)$ is QAR(1) as defined in~(\protect\ref{eqnmodel1}); in the right-hand
picture, $Y_t = X_t / \Var(X_t)^{1/2}$ where $(X_t)$ is the ARCH(1)
process defined in~(\protect\ref{eqnmodel3}). All curves are plotted against
$\omega/ (2\uppi)$.}\label{figL2Spectra}
\end{figure}

%
%
\begin{figure}

\includegraphics{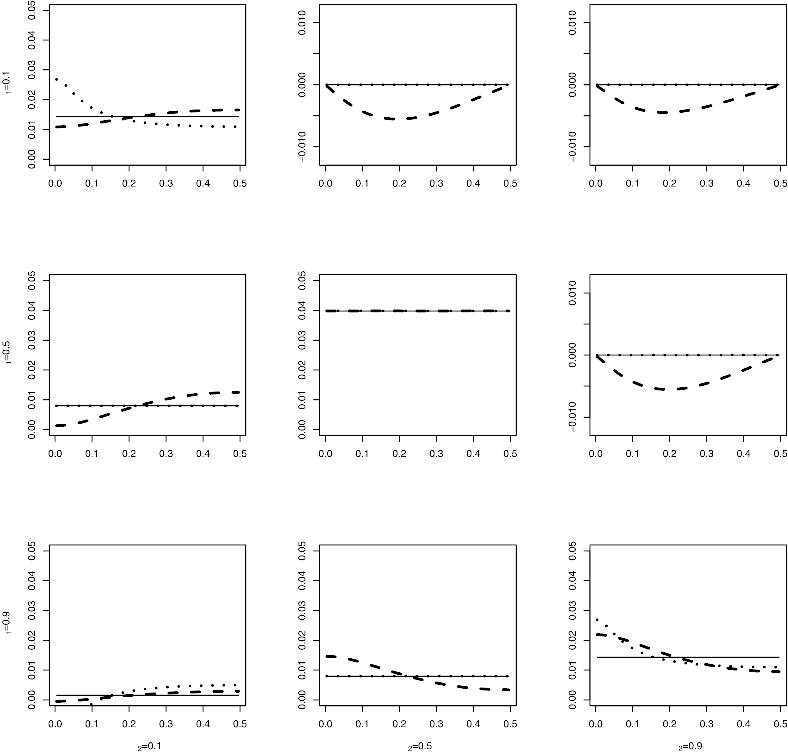}

\caption{Copula spectra ${(2\uppi)^{-1} \sum_{k \in\IZ} \Cov(I\{
F(Y_{t+k}) \leq\tau_1\}, I\{F(Y_t) \leq\tau_2\}) \mathrm
{e}^{-\mathrm{i} \omega k}}$ for $\tau_1,\tau_2 = 0.1, 0.5$, and 0.9.
Real parts (imaginary parts)
are shown in sub-figures with~$\tau_2 \leq\tau_1$ ($\tau_2 > \tau
_1$). Solid, dashed, and dotted lines correspond to the white noise,
QAR(1) and ARCH(1) processes in Figure~\protect\ref{figL2Spectra}.
All curves
are plotted against $\omega/ (2\uppi)$.}\label{figCopulaSpectra}
\end{figure}

As an illustration, the $L^2$-spectra and copula spectral densities are
shown in Figures~\ref{figL2Spectra} and~\ref{figCopulaSpectra},
respectively, for three different processes: (a) a Gaussian white noise process,
(b) a QAR(1) process~(Koenker and Xiao \cite
{Koenker2006}) and (c) an ARCH(1) process
[the same processes are also considered in the simulations of
Section~\ref{secsim}]. All processes were standardized so that the
marginal distributions have unit variance. Although their dynamics
obviously are quite different, those three processes are uncorrelated,
and thus all exhibit the same flat $L^2$-spectrum.
This very clearly appears in Figure~\ref{figL2Spectra}. In
Figure~\ref{figCopulaSpectra}, the copula spectral densities
associated with various values of $\tau_1$ and $\tau_2$ are shown for
the same processes.
Obviously, the three copula spectral densities differ considerably from
each other and, therefore, provide a much richer information about the
dynamics of those three processes.

For an interpretation of Figure~\ref{figCopulaSpectra}, recall~(\ref
{fqq1}) and~(\ref{fqq2}), and note that
\begin{eqnarray*}
-\gamma_k^U(\tau_1, \tau_2) & =&
- \Cov\bigl(I\{U_t \leq\tau_1\}, I\{ U_{t-k}
\leq\tau_2\}\bigr)
\\
& =& \Cov\bigl(I\{U_t \leq\tau_1\}, I
\{U_{t-k} > \tau_2\}\bigr)
\\
& =& \IP(X_t \leq q_{\tau_1}, X_{t-k} >
q_{\tau_2}) - \tau_1 (1 - \tau_2)
\\
& =& \IP(X_t > q_{\tau_1}, X_{t-k} \leq
q_{\tau_2}) - (1-\tau_1) \tau_2.
\end{eqnarray*}
Hence, $\gamma_k^U(\tau_1, \tau_2)$ is the probability for $\{X_t\}$
to switch from the upper $\tau_2$ tail to the lower $\tau_1$ tail in
$k$ steps, minus the corresponding probability for white noise, which
is also the probability for $\{X_t\}$ to switch from the lower $\tau
_2$ tail to the upper $\tau_1$ tail in $k$ steps, minus the
corresponding probability for white noise.

Copula spectral density kernels, as represented in Figure~\ref{figCopulaSpectra}, thus provide information on those
quantile-crossing, or tail-switching probabilities. In particular, a
non-vanishing imaginary part for $\mathfrak{f}_{q_{\tau_1}, q_{\tau
_2}}(\omega)$ indicates that $\IP(X_t \leq q_{\tau_1}, X_{t-k} \leq
q_{\tau_2})$, for some values of~$k$, differs from $\IP(X_t \leq
q_{\tau_1}, X_{t+k} \leq q_{\tau_2})$, which implies that $\{X_t\}$
is not time-revertible. Figure~\ref{figCopulaSpectra}, where
imaginary parts are depicted above the diagonal, clearly indicates that
the QAR(1) process is not time-revertible.

Note that, in order to distinguish between the ARCH(1) and the i.i.d.
process, it is common practice to compare the $L_2$-spectral densities
of the squared processes. This approach can also be used for the QAR(1)
process, but is bound to miss important features. For example, the
asymmetric nature of QAR(1) dynamics, revealed, for example, by the
difference between its $(0.1, 0.1)$ (top left panel) and $(0.9, 0.9)$
(bottom right) spectra cannot be detected in the $L^2$-spectrum of a
squared QAR(1) process.

For a more detailed discussion of the advantages of the copula spectrum
compared to the classical one, see Hong \cite
{Hong2000}, Dette \textit{et~al.} \cite{DetteEtAl2013},
Hagemann \cite{Hagemann2013} and Lee and Rao \cite{LeeRao2012}.

Consistent estimation of $\mathfrak{f}_{q_{\tau_1}, q_{\tau
_2}}(\omega)$ was independently considered in Hagemann \cite{Hagemann2013}
for the special case $\tau_1=\tau_2\in(0,1)$, and by Dette \textit
{et~al.} \cite{DetteEtAl2013} for general couples $(\tau
_1,\tau_2) \in(0,1)^2$ of
quantile levels, under different assumptions such as $m(n)$-decomposa\-
bility and $\beta$-mixing.

Hagemann's estimator, called the \textit{$\tau$th quantile periodogram},
is a traditional $L^2$-peri\-odogram where observations are replaced
with the indicators
\[
I\bigl\{\hat F_n(X_t) \leq\tau\bigr\} = I\{
R_{n;t} \leq n\tau\},
\]
where $\hat F_n(x):= n^{-1}\sum_{t=0}^{n-1} I\{X_t\leq x\}$ denotes
the empirical marginal distribution function and $R_{n;t}$ the rank
of~$X_t$ among $X_0,\ldots, X_{n-1}$. Dette \textit{et~al.}
\cite{DetteEtAl2013} introduce
their Laplace rank-based periodograms by substituting an $L^1$-approach
for the $L^2$ one, and considering the cross-periodograms associated
with arbitrary couples~$(\tau_1,\tau_2)$ of quantile levels. See
Remark~\ref{R21} for details.

In this paper, we stick to the $L^2$-approach, but extend Hagemann's
concept by considering, as in Dette \textit{et~al.} \cite{DetteEtAl2013}, the
cross-periodograms associated with arbitrary couples $(\tau_1,\tau
_2)$. More precisely, we define the \textit{rank-based copula
periodogram}~$I_{n,R}$, shortly, the \textit{CR-periodogram} as the collection
\begin{equation}
\label{eqinr} \inr^{\tau_1, \tau_2}(\omega):= \frac{1}{2\uppi n}
d_{n,R}^{\tau_1}(\omega) d_{n,R}^{\tau_2}(-
\omega), \qquad\omega \in\IR, (\tau_1,\tau_2)
\in[0,1]^2,
\end{equation}
with
\[
d_{n,R}^{\tau}(\omega):= \sum_{t=0}^{n-1}
I\bigl\{\hat F_n(X_t) \leq \tau\bigr\}
\ee^{- \ii\omega t} = \sum_{t=0}^{n-1} I\{
R_{n;t} \leq n\tau\} \ee^{- \ii\omega t}.
\]
Those cross-periodograms, as well as Hagemann's $\tau$th quantile
periodograms, are measurable functions of the marginal ranks $R_{n;t}$,
whence the terminology and the notation.

Classical periodograms and rank-based Laplace periodograms converge, as
$n\to\infty$, to random variables whose expectations are the
corresponding spectral densities; but they fail estimating those
spectral densities in a consistent way.
Similarly, the CR-periodogram $\inr^{\tau_1, \tau_2}(\omega)$ fails
to estimate~$\mathfrak{f}_{q_{\tau_1}, q_{\tau_2}}(\omega)$
consistently. More precisely, let $ \weak$ denote the \textit
{Hoffman--J{\o}rgensen convergence}, namely, the weak convergence
in the space of bounded functions $[0,1]^2\to\mathbb{C}$, which we denote by $\ell_{\mathbb{C}}^\infty
([0,1]^2)$.
Note that results in empirical process theory are typically stated for spaces of
real-valued, bounded functions; see Chapter 1 of van der Vaart and Wellner \cite{vanderVaartWellner1996}.
By identifying $\ell_{\mathbb{C}}^\infty ([0,1]^2)$
with the product space $\ell^\infty ([0,1]^2)\times \ell^\infty ([0,1]^2)$
these results transfer immediately.
We show (see Proposition~\ref{propinr} for details) that, under
suitable assumptions, for any fixed frequencies~$\omega\neq0 \modt
2\uppi$,
\[
\bigl( \inr^{\tau_1, \tau_2}(\omega) \bigr)_{ (\tau_1,\tau_2) \in
[0,1]^2} \weak \bigl(\PerProc(
\tau_1,\tau_2;\omega) \bigr)_{ (\tau_1,\tau_2) \in[0,1]^2} \qquad
\mbox{as } n\to\infty,
\]
where the limiting process $\PerProc$ is such that
\[
\IE\bigl[\PerProc(\tau_1,\tau_2;\omega)\bigr] =
\mathfrak{f}_{q_{\tau_1},
q_{\tau_2}}(\omega) \qquad\mbox{for all }(\tau_1,
\tau_2) \in [0,1]^2\mbox{ and }\omega\neq0 \mod2\uppi
\]
and $\PerProc(\tau_1,\tau_2;\omega_1)$ and $\PerProc(\tau_3,\tau
_4;\omega_2)$ are independent for any $\tau_1,\ldots,\tau_4$ as
soon as $\omega_1\neq\omega_2$.

In view of this asymptotic independence
at different frequencies, it seems natural to consider smoothed
versions of $\inr^{\tau_1, \tau_2}(\omega)$, namely, for $ (\tau
_1, \tau_2) \in[0,1]^2$ and $ \omega\in\mathbb{R} $, averages of
the form
\begin{equation}
\label{eqnDefRankEstimator} \hat G_{n,R}(\tau_1, \tau_2;
\omega):= \frac{2\uppi}{n} \sum_{s=1}^{n-1}
W_n ( \omega- 2\uppi s / n ) I_{n,R}^{\tau_1,
\tau_2}(2\uppi s /
n),
\end{equation}
where $W_n$ denotes a sequence of weighting functions. For the special
case $\tau_1=\tau_2$, the consistency of a closely related estimator
is established by Hagemann \cite{Hagemann2013}.
However, even for $\tau_1 =
\tau_2$, obtaining the asymptotic distribution of \textit{smoothed
CR-periodograms} of the form~(\ref{eqnDefRankEstimator}) is not
trivial, and so far has remained an open problem. Similarly,
Dette \textit{et~al.} \cite{DetteEtAl2013} do not provide
any results on the asymptotic
distributions of their (smoothed) \textit{Laplace rank-based
periodograms}. Note that even consistency results in Hagemann \cite
{Hagemann2013}, as well as in Dette \textit
{et~al.} \cite{DetteEtAl2013} are only pointwise
in~$\tau_1,\tau_2$.

In the present paper, we fill that gap. Theorem~\ref
{thmAsympDensityRankEstimator} below does not only provide pointwise
asymptotic distributions for smoothed CR-periodograms, but also describes
the asymptotic \mbox{behavior} of a properly centered and rescaled version of
the full collection~$\{\hat G_{n,R}(\tau_1, \tau_2; \omega),   (\tau
_1,\tau_2) \in[0,1]^2\}$
as a sequence of \textit{stochastic processes}. Such convergence
results (process convergence rather than pointwise) are of particular
importance, as they can be used to obtain the asymptotic distribution
of functionals of smoothed CR-periodograms as estimators of functionals
of the corresponding copula spectral density kernel. As an example, we
derive, in Section~\ref{secrankper}, the asymptotic distributions of
periodograms computed from various rank-based autocorrelation concepts
(Spearman, Gini, Blomqvist, etc.).

In the process of analyzing the asymptotic behavior of $\{\hat
G_{n,R}(\tau_1, \tau_2; \omega)\}$, we establish several
intermediate results of independent interest. For instance, we prove an
asymptotic representation theorem (Theorem~\ref
{thmAsympDensityEstimator}(i)), where we show that, uniformly in $\tau
_1, \tau_2 \in[0,1]^2,  \omega\in\mathbb{R}$, the smoothed
periodogram~$\hat G_{n,R}(\tau_1, \tau_2; \omega)$ can be approximated~by\vspace*{-1pt}
\begin{equation}
\hat G_{n,U}(\tau_1, \tau_2; \omega):=
\frac{2\uppi}{n} \sum_{s=1}^{n-1}
W_n ( \omega- 2\uppi s / n ) I_{n,U}^{\tau_1,
\tau_2}(2\uppi s /
n),
\end{equation}
where\vspace*{-1pt}
\begin{equation}
\label{eqdefdnu} I_{n,U}^{\tau_1, \tau_2}(\omega):= \frac{1}{2\uppi}
\frac{1}{n} d_{n,U}^{\tau_1}(\omega) d_{n,U}^{\tau_2}(-
\omega),
\end{equation}
and\vspace*{-1pt}
\[
d_{n,U}^{\tau}(\omega):= \sum_{t=0}^{n-1}
I\{U_t \leq\tau\} \ee ^{- \ii\omega t}\qquad\mbox{with }
U_t:= F(X_t).
\]

We conclude this section with two remarks clarifying the relation
between the approach considered here that of Dette \textit
{et~al.} \cite{DetteEtAl2013},
and some other copula-based approaches in the analysis of time series.

%
%
\begin{bem}\label{R21}
The classical $L^2$-periodo\-gram of a real-valued time series can be
represented in two distinct ways, providing two distinct
interpretations. First, it can be defined as the Fourier transform of
the empirical autocovariance function. More precisely, considering the
empirical autocovariance\vspace*{-1pt}
\[
\hat\gamma_k:= \frac{1}{n-k} \sum_{t=1}^{n-k}(X_{t+k}-
\bar X) (X_t - \bar X), \qquad k\geq0, \hat\gamma_k:= \hat\gamma_{-k}, k<0,
\]
the classical $L^2$-periodogram can be defined as
\begin{equation}
\label{defperL2v1} I_n(\omega):= \frac{1}{2\uppi} \sum
_{\llvert  k\rrvert   < n} \frac{n-k}{n} \hat \gamma_k
\ee^{-\ii k \omega}.
\end{equation}
However, an alternative definition
is
\begin{equation}
\label{defperL2v2} I_n(\omega):= \frac{1}{2\uppi} \frac{1}{n}
\Biggl\llvert \sum_{t=0}^{n-1}X_t
\ee^{-\ii t \omega}\Biggr\rrvert ^2 = \frac{n}{4}\bigl(\hat
b_1^2 + \hat b_2^2\bigr),
\end{equation}
where $b_1,b_2$ are the coefficients of the $L^2$-projection of the
observations $X_0,\ldots,X_{n-1}$ on the basis $(1,\sin(\omega
t),\cos
(\omega t))$, that is,
\begin{equation}
\label{defperL2v2loss} (\hat a, \hat b_1,\hat b_2) =
\mathop{\operatorname{Argmin}}_{(a,b_1,b_2)'\in
\mathbb
{R}^3} \sum_{t=0}^{n-1}
\bigl(X_t - a -b_1\cos(\omega t) - b_2\sin
(\omega t) \bigr)^2.
\end{equation}

This suggests two different starting points for generalization.
We either can replace autocovariances in (\ref{defperL2v1}) by
alternative measures of dependence such as (empirical) joint
distributions or copulas, or consider alternative loss functions in the
minimization step~(\ref{defperL2v2loss}). Replacing the
autocovariance function by the pairwise copula with $\tau_1=\tau_2
=\tau$ yields the $\tau$-quantile periodogram of Hagemann \cite
{Hagemann2013}, which we also consider here,
under the name of
CR-periodogram, albeit for general $(\tau_1,\tau_2) \in[0,1]^2$.
Replacing the quadratic loss in~(\ref{defperL2v2loss}) was, in a time
series context, first considered by Li
\cite{Li2008,Li2012} and Dette \textit{et~al.} \cite
{DetteEtAl2013}, who\vspace*{1pt} observed that substituting the check
function~$\rho_\tau(x) = x(\tau- I\{x<0\})$ of Koenker and
Bassett \cite{KoenkerBassett1978} for the standard
$L^2$-loss leads to an estimator
for the quantity
\[
\tilde{\mathfrak{f}}_{\tau, \tau}(\omega):= \frac{1}{2\uppi f^2(q_{\tau})} \sum
_{k\in\IZ} \mathrm{e}^{-\ii
\omega
k} \bigl(\IP(X_0
\leq q_{\tau},X_{ -k } \leq q_{\tau}) -
\tau^2 \bigr).
\]
This latter expression is a weighted version of the copula spectral
density kernel at $\tau_1 = \tau_2 = \tau$ introduced in (\ref
{fqq2}). This weighting, which involves~$ f(q_{\tau})$, is undesirable,
since it involves the unknown marginal distribution of $X_t$, which is
unrelated with its dynamics. Dette \textit{et~al.} \cite
{DetteEtAl2013} demonstrate that,
by considering ranks instead of the original data, that weighting can
be removed. The same authors also proposed a generalization to
cross-periodograms associated with distinct quantile levels. See
Li \cite{Li2012}, Dette \textit{et~al.}
\cite{DetteEtAl2013} and Hagemann \cite
{Hagemann2013} for details
and discussion.
\end{bem}

%
%
\begin{bem}
{The benefits of considering joint distributions and copulas as measures
of serial dependence in a nonparametric time-domain analysis of time
series has been realized by many authors. Skaug and
Tj{\o}stheim \cite{SkaugTjostheim1993}
and Hong \cite{Hong1999} used joint distribution
functions to test for
serial independence at given lag. Subsequently, related approaches were
taken by many authors, and an overview of related results can be found
in Tj{\o}stheim \cite{Tjostheim1996} and
Hong \cite{Hong1999}. Copula-based tests of
serial independence were considered by Genest and
R{\'e}millard \cite{GenestRemillard2004},
among others. Linton and Whang \cite{LintonWhang2007}
introduced the so-called \textit
{quantilogram}, defined as the autocorrelation of the series of
indicators~$I\{X_t \leq\hat q_\tau\}$, $t=0,\ldots, n-1$, where
$\hat q_\tau$ denotes the empirical~$\tau$-quantile; they discuss the
application of this quantilogram (closely related to Hagemann's $\tau
$-quantile periodogram) to measuring directional predictability of time
series. They do not, however, enter into any spectral considerations.
An extension of those concepts to the dependence between several time
series was recently considered in Han \textit{et~al.}
\cite{HanLintonOkaWhang2013}.
Finally, Davis and Mikosch \cite{DavisMikosch2009}
also considered a related quantity
which is based on autocorrelations of indicators of extreme events.
}
\end{bem}

\section{Asymptotic properties of rank-based copula periodograms}\label{secasy}

The derivation of the asymptotic properties of CR-periodograms requires
some assumptions on the underlying process and the weighting functions~$W_n$.

Recall that the \textit{$r$th order joint cumulant} $\cum(\zeta_1,
\ldots, \zeta_r)$ of the random vector $(\zeta_1, \ldots, \zeta
_r)$ is defined as
\[
\cum(\zeta_1, \ldots, \zeta_r):= \sum
_{\{\nu_1, \ldots, \nu_p\}
} (-1)^{p-1} (p-1)! \biggl(\IE\prod
_{j \in\nu_1} \zeta_j\biggr) \cdots\biggl(\IE \prod
_{j \in\nu_p} \zeta_j\biggr),
\]
with summation extending over all partitions $\{\nu_1, \ldots, \nu
_p\}$, $p=1,\ldots,r$ of~$\{1,\ldots,r\}$ (cf. Brillinger \cite
{Brillinger1975}, page~19).

The assumption we make on the dependence structure of the
process~$(X_t)_{t\in\IZ}$ is as follows. Its relation to more
classical assumptions of weak dependence is discussed in
Propositions~\ref{propalph} and~\ref{propgmc} below, and in
Lemma~\ref{lemgmc}.
\begin{itemize}[(C)]
\item[(C)] \textit{There exist constants $\rho\in(0,1)$ and $K <
\infty$ such that}, \textit{for arbitrary intervals $A_1,\ldots, A_p \subset
\mathbb{R}$ and arbitrary $t_1,\ldots,t_p \in\IZ$},
\end{itemize}
%
\begin{equation}
\label{eqboundcum} \bigl\llvert \cum\bigl(I\{X_{t_1} \in A_1\},
\ldots, I\{X_{t_p} \in A_p\}\bigr)\bigr\rrvert \leq K
\rho^{\max_{i,j} \llvert  t_i-t_j\rrvert  }.
\end{equation}

The crucial point here is that we replace an assumption on the
cumulants of the original observations by an assumption on the
cumulants of certain indicators. Thus, in contrast to classical
assumptions, condition \textup{(C)} does not require the existence of any
moments. Additionally, note that the sets~$A_j$ in~(\ref{eqboundcum})
only need to be intervals, not general Borel sets as in classical
mixing assumptions.

%
%
\begin{prop}\label{propalph}
Assume that the process $(X_t)_{t\in\IZ}$ is strictly stationary and
exponentially \mbox{$\alpha$-}mixing, that is,
\begin{equation}
\label{eqalph} \alpha(n):= \mathop{\sup_{A \in\sigma(X_0,X_{-1},\ldots)}}_{B
\in\sigma(X_{n}, X_{n+1},\ldots)}
\bigl\llvert \IP(A \cap B) - \IP(A)\IP(B) \bigr\rrvert \leq K\kappa^n,\qquad n\in\IN
\end{equation}
for some $K<\infty$ and $\kappa\in(0,1)$. Then assumption~\textup{(C)} holds.
\end{prop}

While mixing assumptions are very general and intuitively
interpretable, which makes them quite attractive from a probabilistic
point of view, verifying conditions such as (\ref{eqboundcum}) or
(\ref{eqalph}) can be difficult in specific applications. An
alternative description of dependence that is often easier to check for
was recently proposed by Wu and Shao \cite{wush2004}.
More precisely, these
authors assume that the process~$(X_t)_{t\in\IZ}$ can be represented as
\begin{equation}
\label{defGMC} X_t = g(\ldots,\eps_{t-2},
\eps_{t-1},\eps_{t}), \qquad t \in\IN,
\end{equation}
where $g$ denotes some measurable function and $(\eps_t)_{t \in\IZ}$
is a collection of i.i.d. random variables. Note that the function $g$
is not assumed to be linear, which makes this kind of process very general.
To quantify the long-run dependence between $(\ldots,X_{-1},X_0)$, and
$(X_t,X_{t+1},\ldots)$, denote by $(\eps_t^*)_{t \leq0}$ an independent
copy of $(\eps_t)_{t \leq0}$ and define
\[
X_t^*:= g\bigl(\ldots,\eps_{-1}^*,\eps_0^*,
\eps_1,\ldots,\eps_t\bigr), \qquad t \in\IN.
\]
The process $(X_t)_{t \in\IZ}$ satisfies a \textit{geometric moment
contraction of order $a$} property (shortly, GMC(a) throughout this
paper) if, for some~$K<\infty$ and~$ \sigma\in(0,1)$,
\begin{equation}
\label{eqgmc} \E\bigl\llvert X_n - X_n^*\bigr\rrvert
^a \leq K \sigma^n,\qquad n \in\IN;
\end{equation}
see Wu and Shao \cite{wush2004}.
Examples of processes that satisfy this condition include, (possibly,
under mild additional conditions on the parameters) ARMA, ARCH, GARCH,
asymmetric GARCH, random coefficient autoregressive, quantile
autoregressive and Markov models, to name just a few. Proofs and
additional examples can be found in Shao and Wu \cite
{ShaoWu2007} and Shao \cite{Shao2010WN}. The
definition in (\ref{eqgmc}) still requires the
existence of moments, which is quite undesirable in our setting.
However, the following result shows that a modified version of (\ref
{eqgmc}) is sufficient for our purposes.

%
%
\begin{prop}\label{propgmc}
Assume that the strictly stationary process $(X_t)_{t\in\IZ}$ can be
represented as in~(\ref{defGMC}), and that $X_0$ has distribution
function $ F$. Let the process $(F(X_t))_{t \in\IZ}$ satisfy GMC(a)
for some $a>0$, that is, assume that there exist $K<\infty$ and $\sigma
\in(0,1)$ such that
\begin{equation}
\label{eqfgmc} \E\bigl\llvert F(X_n) - F\bigl(X_n^*
\bigr)\bigr\rrvert ^a \leq K \sigma^n,\qquad n \in\IN.
\end{equation}
Then assumption~\textup{(C)} holds.
\end{prop}

The important difference between assumptions (\ref{eqgmc}) and (\ref
{eqfgmc}) lies in the fact that, in condition~(\ref{eqfgmc}), only
the random variables $F(X_t)=U_t$, which possess moments of arbitrary
order, appear. This implies that a GMC(a) condition on $X_t$ with
arbitrarily small values of $a$, together with a very mild regularity
condition on $F$, are sufficient to imply assumption~\textup{(C)}. More
precisely, we have the following result.

%
%
\begin{lemma}\label{lemgmc}
Assume that $(X_t)_{t\in\IZ}$ is strictly stationary.
Let $(X_t)_{t\in\IZ}$ satisfy the GMC(b) condition for some $b >
0$, and assume that the distribution function $F$ of $X_0$ is H{\"o}lder-continuous of order $\gamma>0$. Then (\ref{eqfgmc}) holds for
any $a>0$.
\end{lemma}

For a proof of Lemma \ref{lemgmc}, note that
\[
\IE\bigl\llvert F(X_t)-F\bigl(X_t^*\bigr) \bigr
\rrvert ^a \leq2^{a-b/\gamma}\IE \bigl\llvert F(X_t)-F
\bigl(X_t^*\bigr) \bigr\rrvert ^{b/\gamma} \leq C \IE\bigl
\llvert X_t - X_t^*\bigr\rrvert ^{b} \leq C
K \sigma^{t},
\]
where $\sigma\in(0,1)$ and $K > 0$ are the constants from the GMC(b)
condition.

%
%
\begin{bem}
Although not very deep at first sight, the above result has some
remarkable implications. In particular, we show in the \hyperref[pmt]{Appendix} that,
under a very mild regularity condition on~$F$, the copula spectra of a
GMC(a) process are analytical functions of the frequency $\omega$.
This is in sharp contrast with classical spectral density analysis,
where higher-order moments are required to obtain smoothness of the
spectral density.
\end{bem}

We now are ready to state a first result on the asymptotic properties
of the CR-periodogram $\inr^{ \tau_1, \tau_2 }$ defined in (\ref{eqinr}).\vspace*{-2pt}

%
%
\begin{prop}\label{propinr}
Assume that $F$ is continuous and that $(X_t)_{t\in\IZ}$ is strictly
stationary and satisfies assumption~\textup{(C)}. Then, for every fixed $\omega
\neq0 \modt2\uppi$,\vspace*{-1pt}
\[
\bigl( \inr^{\tau_1, \tau_2}(\omega) \bigr)_{ (\tau_1,\tau_2) \in
[0,1]^2} \weak \bigl(\PerProc(
\tau_1,\tau_2;\omega) \bigr)_{ (\tau_1,\tau_2) \in[0,1]^2} \qquad
\mbox{in } \ell^\infty_{\mathbb{C}}\bigl([0,1]^2\bigr).
\]
The (complex-valued) limiting processes $\PerProc$, indexed by $(\tau
_1, \tau_2) \in[0,1]^2$, are of the form\vspace*{-1pt}
\[
\PerProc(\tau_1,\tau_2;\omega) = \frac{1}{2\uppi}
\mathbb{D}(\tau _1;\omega) \overline{\mathbb{D}(\tau_2;
\omega)}
\]
with
$\mathbb{D}(\tau;\omega) = \mathbb{C}(\tau;\omega) + \ii\mathbb
{S}(\tau;\omega)$ where $\mathbb{C}$ and $\mathbb{S}$ denote two
centered jointly Gaussian processes.
For $\omega\in \mathbb{R}$, their covariance structure takes the form\vspace*{-1pt}
\[
\IE \bigl[\bigl(\mathbb{C}(\tau_1;\omega),\mathbb{S}(
\tau_1;\omega )\bigr)^\prime\bigl(\mathbb{C}(
\tau_2;\omega),\mathbb{S}(\tau_2;\omega)\bigr) \bigr] =
\upi\pmatrix{ \displaystyle\Re\mathfrak{f}_{q_{\tau_1},q_{\tau_2}}(\omega) & \displaystyle-
\Im\mathfrak {f}_{q_{\tau_1},q_{\tau_2}}(\omega)
\vspace*{3pt}\cr
\displaystyle\Im
\mathfrak{f}_{q_{\tau_1},q_{\tau_2}}(\omega) & \displaystyle\Re\mathfrak
{f}_{q_{\tau_1},q_{\tau_2}}(\omega)}.
\]
Moreover, $\mathbb{D}(\tau;\omega)=\mathbb{D}(\tau;\omega+2\uppi)=\overline{\mathbb{D}(\tau;-\omega)}$,
and the family $\{\mathbb{D}(\cdot;\omega):\omega\in [0,\uppi]\}$
is a collection of independent processes.\vspace*{-1pt}
\end{prop}

Note\vspace*{1pt} that, for $\omega= 0 \modt2\uppi$ we have $d_{n,R}^{\tau
}(0) = n\tau+ \mathrm{o}_P(n^{1/2})$, where the exact form of the remainder
term depends on the number of ties in the observations. Therefore,
under the assumptions of Proposition~\ref{propinr}, $\inr^{\tau_1,
\tau_2}(0) = n(2\uppi)^{-1}\tau_1\tau_2 + \mathrm{o}_P(1)$ for $\omega= 0
\modt2\uppi$.

%
%
\begin{bem}
Proposition~\ref{propinr} implies that the CR-periodograms
corresponding to different frequencies are asymptotically independent.
Therefore, it can be used to obtain the asymptotic distribution of the
smoothed periodogram for any $b_n = m / n$, with $m \in\IN$ \textit
{not depending on} $n$. In this case, $m$ CR-periodograms are used for
smoothing and the asymptotic distributions of the smoothed
CR-periodograms follow from Proposition~\ref{propinr}. However, in
this case, the smoothed periodogram is not necessarily a consistent
estimator (its variance does not tend to zero) of the spectral density
kernel $\mathfrak{f}_{q_{\tau_1}, q_{\tau_2}}(\omega)$. As we shall
see, in order to be consistent, a smoothed CR-periodogram requires $n b_n
= m(n)  \rightarrow\infty$ as $n \rightarrow\infty$.\vspace*{-1pt}
\end{bem}

In order to establish the convergence of the smoothed CR-periodogram
process (\ref{eqnDefRankEstimator}), we require the weights $W_n$ in
(\ref{eqnDefRankEstimator}) to satisfy the following assumption,
which is quite standard in classical time series analysis (see, e.g.,
page~147 of Brillinger \cite{Brillinger1975}).
\begin{itemize}[(W)]
\item[(W)] The weight function $W$ is real-valued and even, with
support $[-\upi,\upi]$; moreover, it has bounded variation, and satisfies
$
\int_{-\upi}^{\upi} W(u) \,\mathrm{d}u = 1$.
\end{itemize}

Denoting by $b_n > 0$, $n=1,2,\ldots,$ a sequence of scaling
parameters such that~$b_n \rightarrow0$ and $n b_n \rightarrow\infty
$ as $n \rightarrow\infty$, define\vspace*{-4pt}
\[
W_n(u):= \sum_{j=-\infty}^{\infty}
b_n^{-1} W\bigl(b_n^{-1} [u + 2\uppi j]
\bigr).
\]
We now are ready to state our main result.

%
%
\begin{theorem} \label{thmAsympDensityRankEstimator}
Let assumptions \textup{(C)} and \textup{(W)} hold. Assume that $X_0$ has a continuous
distribution function $F$ and that there exist constants $\kappa> 0$
and $k \in\IN$, such that
\[
b_n = \mathrm{o}\bigl(n^{-1/(2k+1)}\bigr)\quad\mbox{and}\quad
b_n n^{1-\kappa} \rightarrow\infty.
\]
Then, for any fixed $\omega\in\IR$, the process
\[
\mathbb{G}_{n}(\cdot,\cdot;\omega):= \sqrt{n b_n}
\bigl(\hat G_{n,R}( \tau_1, \tau_2; \omega) -
\mathfrak{f}_{q_{\tau_1},
q_{\tau_2}}(\omega) - B_n^{(k)}(
\tau_1, \tau_2; \omega) \bigr)_{\tau_1,\tau_2 \in[0,1]}
\]
satisfies
\begin{equation}
\label{eqgntoh} \mathbb{G}_{n}(\cdot,\cdot;\omega) \rightsquigarrow
H(\cdot,\cdot;\omega)
\end{equation}
in $\ell^{\infty}_{\mathbb{C}}([0,1]^2)$, where the bias $B_n^{(k)}$ is given by
\begin{equation}
\label{defbias} B_n^{(k)}(\tau_1,
\tau_2; \omega):=
\sum_{j=2}^k \frac{b_n^j}{j!} \int
_{ -\upi}^{\upi} v^j W(v) \,\dd v
\frac
{\mathrm{d}^j}{\mathrm{d}\omega^j}\mathfrak{f}_{q_{\tau_1},q_{
\tau
_2}}(\omega), \qquad \omega\in \mathbb{R},
\end{equation}
and\vspace*{1pt} $\mathfrak{f}_{q_{\tau_1},q_{ \tau_2}}$ is defined in (\ref
{fqq2}). The process $H(\cdot,\cdot;\omega)$ in~(\ref{eqgntoh}) is
a centered Gaussian process characterized by
\begin{eqnarray*}
&& \Cov \bigl(H(u_1, v_1; \omega ), H( u _2, v
_2; \omega)\bigr)
\\
&&\quad = 2\uppi \biggl(\int_{-\upi}^\upi
W^2(w)\,\mathrm{d}w \biggr)
\bigl(\mathfrak{f}_{q_{u_1}, q_{u_2}}(\omega) \mathfrak {f}_{q_{v_2}, q_{v_1}}(
\omega) + \mathfrak{f}_{q_{u_1}, q_{v_2}}(\omega) \mathfrak{f}_{q_{v_1},
q_{u_2}}(
\omega) I\{\omega=0 \mod\upi\} \bigr).
\end{eqnarray*}
Moreover, $H(\omega)= H(\omega+2\uppi)= \overline{H(-\omega)}$, and the family~$\{H(\omega),  \omega\in[0,\upi] \}
$ is a collection of independent processes. In particular, the weak
convergence~(\ref{eqgntoh}) holds jointly for any finite fixed
collection of frequencies~$\omega$.
\end{theorem}

%
%
\begin{bem}
Assume that $W$ is a kernel of order $d$, that is, $\int_{-\upi}^{\upi
} v^j W(v) \,\mathrm{d}v = 0$, for $j < d$ and $0 < \int_{-\upi}^{\upi} v^d
W(v) \,\mathrm{d}v < \infty$. The Epanechnikov kernel, for example, is of
order 2. Then, for $\omega\neq0 \modt2\uppi$, the bias is of order
$b_n^d$. Since the variance is of order $(n b_n)^{-1}$, the mean
squared error will be minimized when $b_n$ decays at rate
$n^{-1/(2d+1)}$. Therefore, for kernels of finite order, the optimal
bandwidth fulfils the assumptions of Theorem~\ref
{thmAsympDensityRankEstimator}.
\end{bem}

%
%
\begin{bem}\label{bemAsyCI}
Theorem~\ref{thmAsympDensityRankEstimator} can be used to conduct
asymptotic inference in various ways. An important example is the
construction of asymptotic confidence intervals. For illustration
purposes, consider the case $\tau_1 = \tau_2 = \tau$. Assume that
$W$ is a kernel of order $d$, and that the bandwidth~$b_n$ is chosen
such that $b_n = \mathrm{o}(n^{-1/(2d+1)})$. In this case, the bias is of order
$b_n^d$, and thus is asymptotically negligible compared to the
variance. An asymptotic confidence interval thus can be constructed by
using Theorem~\ref{thmAsympDensityRankEstimator} to obtain the approximation
\[
\sqrt{nb_n}\bigl(\hat G_{n,R}( \tau, \tau; \omega) -
\mathfrak {f}_{q_{\tau}, q_{\tau}}(\omega)\bigr) \approx\mathcal{N}\bigl(0,\sigma
^2\bigr)\quad\mbox{for }\sigma^2 = 2\uppi\int
W^2(u) \,\mathrm{d}u \,\mathfrak {f}_{q_{\tau}, q_{\tau}}^2(\omega)
\]
and estimating $\sigma^2$ by plugging in $\hat G_{n,R}( \tau, \tau;
\omega)$ as an estimator for $\mathfrak{f}_{q_{\tau}, q_{\tau
}}(\omega)$. A more detailed discussion of confidence interval
construction that also includes the case where $\tau_1 \neq\tau_2$
is deferred to Section~\ref{secsim}.
\end{bem}

%
%
\begin{bem}
Process convergence with respect to the frequencies $\omega$ cannot
hold since the limiting processes are independent for different values
of $\omega$. This implies that there exists no tight random element in
$\ell^\infty_{\mathbb{C}}([0,1]^2\times[0,\upi])$ with\vspace*{1pt} the right
finite-dimensional distributions, as would be required for process
convergence in $\ell^\infty_{\mathbb{C}}([0,1]^2\times[0,\upi])$ to take place.
Note that a similar situation occurs for the classical $L^2$-spectral
density which does not converge as a process when indexed by frequencies.
\end{bem}

For fixed quantile levels $\tau_1,\tau_2$, the asymptotic
distribution of $\mathbb{G}_{n}(\tau_1,\tau_2;\omega)$ is the same
as the distribution of the smoothed $L^2$-cross-periodogram (see
Chapter~7 of Brillinger \cite{Brillinger1975})
corresponding to the
(unobservable) bivariate time series
\[
\bigl(I\bigl\{F(X_t) \leq\tau_1\bigr\},I\bigl
\{F(X_t) \leq\tau_2\bigr\}\bigr)_{0\leq t\leq n-1}.
\]
In particular, the estimation of the marginal quantiles has no impact
on the asymptotic distribution of $\mathbb{G}_{n}$. Intuitively, this
can be explained by the fact that~$(\hat q_{\tau_1}, \hat q_{\tau
_2})$ converges at~$n^{-1/2}$ rate while the normalization $\sqrt
{nb_n}$ appearing in $\mathbb{G}_{n}$ is strictly slower.

One aspect of Theorem \ref{thmAsympDensityRankEstimator} that does
not appear in the context of classical spectral density estimation is
the convergence of $\mathbb{G}_{n}$ as a process. Establishing this
result is challenging, and it requires the development of new tools. On
the other hand, once convergence has been established at process level,
it can be applied to derive the asymptotic distributions of various
related statistics; see Section~\ref{secrankper}.

%
%
\begin{bem}
In the derivation of Theorem \ref{thmAsympDensityRankEstimator}, it
would be natural to show that~$d_{n,R}^\tau(\omega)$ and
$d_{n,U}^\tau(\omega)$ are sufficiently close to each other uniformly
with respect to $\tau$ and $\omega$, as~$n\to\infty$. Indeed, using
modifications of standard arguments from empirical process theory, it
is possible to establish that
\begin{equation}
\label{approxd} n^{-1/2}\mathop{\sup_{ \omega\in\IR}}_{\tau\in[0,1] }
\bigl\llvert d_{n,R}^\tau(\omega) - d_{n,U}^\tau(
\omega)\bigr\rrvert = \mathrm{o}_P(r_n)
\end{equation}
for some rate $r_n \rightarrow0$ depending on the underlying
dependence structure. Unfortunately, the best rate $r_n$ that can
theoretically be obtained \textit{must} be slower than $\mathrm{o}(n^{-1/4})$,
and this makes the~approximation (\ref{approxd}) useless in
establishing Theorem \ref{thmAsympDensityRankEstimator} for
practically relevant choices of the bandwidth parameter.
\end{bem}

%
%
\begin{bem}
Another type of process convergence is frequently discussed in the
literature on classical $L^2$-based spectral analysis, which is dealing with
\textit{empirical spectral processes} of the form
\[
\biggl( \int_{-\upi}^{\upi} g(\omega) I_n(
\omega)\,\dd\omega \biggr)_{g\in
\mathcal{G}}
\]
with $\mathcal{G}$ denoting a suitable class of functions. For more
details, see Dahlhaus \cite{Dahlhaus1988},
Dahlhaus and Polonik \cite{DahlhausPolonik2009}, and
the references therein. Those processes are completely different from
the processes considered above, and the mathematical tools that need to
be developed for their analysis also differ substantially. It would be
very interesting to extend our results to classes of integrated
periodograms that are indexed by classes of functions. Such an
extension, however, is beyond of the scope of the present paper.
\end{bem}

%
%
\begin{bem}
At first glance, it seems surprising that the asymptotic theory
developed here does not require the marginal distribution function $F$
to have a continuous Lebesgue density, although the CR-periodograms in
(\ref{eqinr}) are based on marginal quantiles. The reason is that the
estimators which are constructed from $X_0,\ldots,X_{n-1}$ are almost
surely equal to estimators based on the (unobserved) transformed
variables $F(X_0),\ldots,F(X_{n-1})$. A similar phenomenon can be observed
in the estimation of copulas; see, for example,
Fermanian, Radulovi{\'c} and Wegkamp
\cite{ferradweg2004}.
\end{bem}

In\vspace*{1pt} order to establish Theorem \ref{thmAsympDensityRankEstimator}, we
prove (an \textit{asymptotic representation} result) that the estimator
$\hat G_{n, R}$ can be approximated by $\hat G_{n, U}$ in a suitable
uniform sense. Theorem~\ref{thmAsympDensityRankEstimator} then
follows from the asymptotic properties of~$\hat G_{n, U}$, which we
state now.

%
%
\begin{theorem}\label{thmAsympDensityEstimator}
Let assumptions~\textup{(C)} and \textup{(W)} hold, and assume that the distribution
function $F$ of $X_0$ is continuous.
Let $b_n$ satisfy the assumptions of Theorem~\ref
{thmAsympDensityRankEstimator}.
Then,
\begin{longlist}[(iii)]
\item[(i)] for any fixed $\omega\in\IR$, as $n\to\infty$,
\[
\sqrt{nb_n} \bigl( \hat G_{n,U}(\tau_1,
\tau_2; \omega) - \IE\hat G_{n,U}(\tau_1,
\tau_2; \omega) \bigr)_{\tau_1, \tau_2 \in[0,1]} \rightsquigarrow H(\cdot,\cdot;
\omega)
\]
in\vspace*{1pt} $\ell^{\infty}_{\mathbb{C}}([0,1]^2)$, where the process $H(\cdot,\cdot;\omega)$ is defined in Theorem \ref{thmAsympDensityRankEstimator};
\item[(ii)] still as $n\to\infty$,
\[
\mathop{\sup_{\tau_1, \tau_2 \in[0,1]}}_{\omega\in\mathbb
{R}} \bigl\llvert \IE\hat
G_{n,U}(\tau_1, \tau_2; \omega) - \mathfrak
{f}_{q_{\tau_1}, q_{\tau_2}}(\omega) - B_n^{(k)}(\tau_1,
\tau_2, \omega)\bigr\rrvert = \mathrm{O}\bigl((nb_n)^{-1}
\bigr) + \mathrm{o}\bigl(b_n^k\bigr),
\]
where $B_n^{(k)}$ is defined in~(\ref{defbias});
\item[(iii)] for any fixed $\omega\in\mathbb{R}$,
\[
\sup_{\tau_1, \tau_2 \in[0,1]} \bigl\llvert \hat G_{n, R} (
\tau_1, \tau_2; \omega) - \hat G_{n, U} (
\tau_1, \tau_2; \omega) \bigr\rrvert = \mathrm{o}_P
\bigl( (n b_n)^{-1/2}+ b_n^k \bigr);
\]
if moreover the kernel $W$ is uniformly Lipschitz-continuous, this
bound is uniform with respect to $\omega\in\mathbb{R}$.
\end{longlist}
\end{theorem}

\section{Spearman, Blomqvist and Gini spectra}\label{secrankper}

In the past decades, considerable effort has been put into replacing
empirical autocovariances by alternative (scalar) measures of
dependence; see, for example, Kendall \cite
{kendall1938}, Blomqvist \cite{blomqvist1950},
Cifarelli, Conti and Regazzini
\cite{CifarelliContiRegazzini1996},
Ferguson, Genest and Hallin
\cite{FergusonEtAl2000} and Schmid \textit{et~al.} \cite{sstgr2010}
for a recent survey. Such
measures of association provide a good compromise between the limited
information contained in autocovariances on one hand, and the fully
nonparametric nature of joint distributions and copulas on the other.

A particularly appealing class of such dependence measures is given by
general rank-based autocorrelations (see Hallin and
Puri \cite{HallinPuri1992,Hallin1994} or
Hallin \cite{Hallin2012} for a survey). The
idea of using ranks in a time-series context is not new. In fact, it is
possible to trace back rank-based measures of serial dependence to the
early developments of rank-based inference: early examples include run
statistics or the serial version of Spearman's
rho (see Wald and Wolfowitz \cite{WaldWolfowitz1943}).
The asymptotics of rank-based autocorrelations are well studied under
the assumption of white noise or, at least, \textit{exchangeability}, and
under contiguous alternatives of serial dependence. An alternative
approach to deriving the asymptotic distribution of rank-based
autocorrelations, which is applicable under general kinds of
dependence, is based on their representation as functionals of
(weighted) empirical copula processes and was considered, for instance,
in Fermanian, Radulovi{\'c} and Wegkamp~\cite{ferradweg2004},
Berghaus, B\"ucher and Volgushev \cite{BuecherEtAl2014}.

Despite the great success of the $L^2$-periodogram in time series
analysis, the only attempt to consider Fourier transforms of rank-based
autocorrelations (or any other rank-based scalar measures of
dependence), to the best of our knowledge, is that of Ahdesm{\"{a}}ki
\textit {et~al.} \cite{Ahdesmaki2005}. The latter paper is of a
more empirical nature, and no
theoretical foundation is provided. The aim of the present section is
to introduce a general class of frequency domain methods, and discuss
their connection to rank-based extensions of autocovariances.

\subsection{The Spearman periodogram} To illustrate our purpose, first
consider in detail the classical example of Spearman's rank
autocorrelation coefficients (more precisely, a version of it~-- see
Remark~\ref{Ahdesmaki}); at lag~$k$, that coefficient can be defined as
\[
\hat\rho_n^k:= \frac{ 12 }{n^3} \sum
_{t=0}^{n-\llvert  k\rrvert  -1} \biggl(R_{n;t} -
\frac{n+1}{2} \biggr) \biggl(R_{n;t+\llvert  k\rrvert  } - \frac{n+1}{2} \biggr).
\]
Letting $ \mathcal{F}_n:= \{ {2\uppi j}/{n}\mid j = 1, \ldots,
\lfloor\frac{n-1}{2}\rfloor- 1, \lfloor\frac{n-1}{2}\rfloor\}$,
define the \textit{Spearman} and \textit{smoothed Spearman
periodograms} as
\[
I_{n,\rho}(\omega):= \frac{1}{2\uppi}\sum_{\llvert  k\rrvert  <n}
\ee^{-\ii\omega k} \hat\rho_n^k, \qquad\omega\in
\mathcal{F}_n
\]
and
\[
\hat G_{n,\rho}(\omega):= \frac{2\uppi}{n} \sum
_{s=1}^{n-1} W_n ( \omega- 2\uppi s / n )
I_{n,\rho}(2\uppi s / n), \qquad\omega \in \mathbb{R},
\]
respectively. Intuition suggests that the (smoothed) rank-based
periodogram~$\hat G_{n,\rho}$ should be an estimator for the Fourier transform
\[
\mathfrak{f}_{\rho}(\omega):= \frac{1}{2\uppi} \frac{1}{12}
\sum_{k\in\IZ} \mathrm{e}^{-\ii
\omega
k}
\rho_k
\]
of the population counterpart
\begin{equation}
\label{defrhok} \rho_k = \rho(C_k)=12\int
_{[0,1]^2} \bigl(C_k(u,v) - uv\bigr) \,\dd u\,\dd v,
\end{equation}
of $\hat\rho_n^k$, where $C_k$ is the copula associated with $(X_t,
X_{t+k})$ (see, e.g., Schmid \textit{et~al.} \cite{sstgr2010}).
Due to the presence of ranks,
the investigation of the asymptotic properties of the Spearman
periodogram under non-exchangeable observations seems highly
non-trivial. However, as we shall demonstrate now, those properties can
be obtained via Theorem~\ref{thmAsympDensityRankEstimator} by
establishing a connection between the Spearman periodogram and the
CR-periodogram.

%
%
\begin{prop}\label{proprho1} For any $\omega\in\mathcal{F}_n $,
\begin{equation}
\label{eqrhofromranks} I_{n,\rho}(\omega) = 12\int_{[0,1]^2}
I_{n,R}^{u,v}(\omega) \,\dd u \,\dd v,
\end{equation}
where $I_{n,R}^{u,v}$ is defined in (\ref{eqinr})
Moreover, for any fixed $\omega\in\mathbb{R}$,
\[
\hat G_{n,\rho}(\omega) = 12\int_{[0,1]^2} \hat
G_{n,R}(u,v;\omega) \,\dd u \,\dd v,
\]
where $\hat G_{n,R}$ is defined in (\ref{eqnDefRankEstimator}).
\end{prop}

\begin{pf*}{Proof of Proposition \ref{proprho1}}
Simple algebra yields
\[
I_{n,\rho}(\omega) = \frac{12}{2\uppi} \frac{1}{n}
d_{n,\rho
}(\omega) d_{n,\rho}(-\omega)\qquad\mbox{with }
d_{n,\rho
}(\omega):= \frac{ 1 }{n} \sum
_{t=0}^{n-1} R_{n;t} \ee^{- \ii
\omega t}.
\]
Observe that
\begin{eqnarray*}
I_{n,\rho}(\omega) =\frac{12}{2\uppi}\frac{1}{n^{ 3}}\sum
_{s,t=0}^{n-1}R_{n;t}R_{n;s}\ee^{-\ii\omega t}\ee^{\ii\omega s}.
\end{eqnarray*}
On the other hand,
\begin{eqnarray}\label{secondterm}
\int_{[0,1]^2}I_{n,R}^{u,v}(\omega)\,\dd u\,\dd v &=&
\frac{12}{2\uppi}\frac
{1}{n}\sum_{s,t=0}^{n-1}
\ee^{-\ii\omega t}\ee^{\ii\omega s} \int_{[0,1]^2} I
\{R_{n;t}\leq n u,R_{n;s}\leq n v\} \,\dd u\,\dd v
\nonumber
\\
&=& \frac{12}{2\uppi}\frac{1}{n}\sum_{s,t=0}^{n-1}
\ee^{-\ii\omega
t}\ee^{\ii\omega s} \bigl(1- n^{-1} R_{n;t}
\bigr) \bigl(1- n^{-1} R_{n;s}\bigr)
\\
&=& I_{n,\rho}(\omega) + \frac{12}{2\uppi}
\frac
{1}{n^2}\sum_{s,t=0}^{n-1}
\ee^{-\ii\omega t}\ee^{\ii\omega s} (n- R_{n;t}- R_{n;s}).\nonumber
\end{eqnarray}
For $\omega\in\mathcal{F}_n$, $\sum_{t=0}^{n-1} \ee^{\ii\omega t} =
0$, so that the second term in (\ref{secondterm}) vanishes. The claim
follows.
\end{pf*}

This result is useful in several ways. On one hand, it allows to easily
derive the asymptotic distribution of the smoothed Spearman periodogram
by applying the continuous mapping theorem in combination with
Theorem~\ref{thmAsympDensityRankEstimator}; see Proposition~\ref
{propAsympSpearmanEstim} below. On the other hand, it motivates the
definition of a general class of rank-based spectra to be discussed in
the next section.

%
%
\begin{prop}\label{propAsympSpearmanEstim}
Under the assumptions of Theorem \ref{thmAsympDensityRankEstimator},
for any fixed frequency \mbox{$\omega\neq0 \modt2\uppi$},
\[
I_{n, \rho}(\omega) \rightsquigarrow12 \int_0^1
\int_0^1 \PerProc (\tau_1,
\tau_2;\omega) \,\dd \tau_1\,\dd\tau_2
\]
and, for every fixed \mbox{$\omega\in\IR$},
\[
\sqrt{n b_n} \bigl(\hat G_{n,\rho}(\omega) -
\mathfrak{f}_{\rho
}(\omega) - B_{n,\rho}^{(k)}(\omega)
\bigr) \stackrel{\mathcal {D}} {\longrightarrow} Z_\rho(\omega),
\]
where $Z_\rho(\omega) \sim\mathcal{N} (0, 2\uppi\frak{f}_{\rho
}^2(\omega) (1+I\{\omega=0 \modt\upi\}) \int W^2(w) \,\mathrm{d}w
)$ and
\[
B_{n,\rho}^{(k)}(\omega):= %
\sum_{j=2}^k \frac{b_n^j}{j!} \int
_{ -\upi}^{\upi} v^j W(v) \,\dd v
\frac
{\mathrm{d}^j}{\mathrm{d}\omega^j}\mathfrak{f}_{\mu}(\omega) , \qquad \omega\in \mathbb{R}.
\]
Moreover, $Z_\rho(\omega) = Z_\rho(-\omega)$, $Z_\rho(\omega) =
Z_\rho(2\uppi+ \omega)$ and $Z_\rho(\omega)$, $\omega\in[0,\upi]$
are mutually independent random variables. The weak convergence above
holds jointly for any finite, fixed collection of frequencies $\omega$.
\end{prop}

This result is a direct consequence of the more general
Proposition~\ref{propAsympGenLin}, which we establish in the next section.
Note that, following the method described in Remark~\ref{bemAsyCI},
Proposition~\ref{propAsympSpearmanEstim} can be used to construct
pointwise asymptotic confidence bands for $\mathfrak{f}_{\rho}(\omega)$.

%
%
\begin{bem}\label{Ahdesmaki}
A closely related version of the Spearman periodogram was recently
considered by Ahdesm{\"{a}}ki \textit
{et~al.} \cite{Ahdesmaki2005}.
The main difference with our
approach is that these authors use a slightly different version of the
lag-$k$ Spearman coefficient, namely
\[
\tilde\rho_k:= \frac{1}{n} \frac{12}{(n-k)^2-1}\sum
_{t=0}^{n-k-1} \biggl(R_{n;t}^k-
\frac{n-k+1}{2} \biggr) \biggl(\bar R_{n;t+k}^k-
\frac{n-k+1}{2} \biggr),
\]
where\vspace*{1pt} $R_{n;t}^k$ denotes the rank of $X_t$ among~$X_0,\ldots,X_{n-k-1}$
and $\bar R_{n;t}^k$ the rank of $X_t$ among $X_{k-1},\ldots,X_{n-1}$,
respectively.\vspace*{1pt} Letting $\tilde\rho_{k}:= \tilde\rho_{-k}$ for
$k<0$, Ahdesm{\"{a}}ki \textit
{et~al.} \cite{Ahdesmaki2005} then
consider a statistic of the form
$\sum_{\llvert  k\rrvert   < n} \ee^{\ii k\omega} \tilde\rho_k$. Note that these
authors investigate their method by means of a simulation study and do
not provide any asymptotic theory.
\end{bem}

\subsection{A general class of rank-based spectra} \label{secgenlin}
The findings in the previous section suggest considering a general
class of rank-based periodograms which are defined in terms of the
CR-periodogram as
\begin{equation}
\label{eqinmu} I_{n,\mu}(\omega):= \int_{[0,1]^2}
I_{n,R}^{u,v}(\omega) \,\mathrm{d}\mu(u,v), \qquad\omega\in
\mathcal{F}_n,
\end{equation}
where $\mu$ denotes an arbitrary finite measure on $[0,1]^2$. A
smoothed version of $I_{n,\mu}$ is defined through
\[
\hat G_{n,\mu}(\omega):= \frac{2\uppi}{n} \sum
_{s=1}^{n-1} W_n ( \omega- 2\uppi s / n )
I_{n,\mu}(2\uppi s / n), \qquad\omega\in \mathbb{R}.
\]
As discussed in the previous section, taking $\mu$ as $12$ times the
uniform distribution on~$[0,1]^2$ yields the Fourier transform of
Spearman autocorrelations.

The general results in Theorem~\ref{thmAsympDensityRankEstimator}
combined with the continuous mapping theorem imply that the smoothed
periodogram $\hat G_{n,\mu}$ is a consistent and asymptotically normal
estimator of a spectrum of the form
\[
\mathfrak{f}_{\mu}(\omega):= \frac{1}{2\uppi} \sum
_{k\in\IZ} \mathrm{e}^{-\ii\omega k} \int_{[0,1]^2}
\bigl(C_k(u,v) - uv\bigr) \,\mathrm{d}\mu(u,v),
\]
where $C_k$ denotes the copula of the pair $(X_0,X_k)$.

%
%
\begin{prop}\label{propAsympGenLin}
Under the assumptions of Theorem~\ref{thmAsympDensityRankEstimator},
for any fixed frequency~\mbox{$\omega\in\IR$},
\[
\sqrt{n b_n} \bigl(\hat G_{n,\mu}(\omega) -
\mathfrak{f}_{\mu
}(\omega) - B_{n,\mu}^{(k)}(\omega)
\bigr) \stackrel{\mathcal {D}} {\longrightarrow} Z_\mu(\omega) \sim
\mathcal{N} \bigl(0, \sigma _\mu^2 \bigr),
\]
where the variance $\sigma^2_\mu$ takes the form
\begin{eqnarray*}
\sigma^2_\mu&=& 2\uppi\int_{-\upi}^\upi
W^2(w) \,\dd w
\\
&&\hspace*{29pt}{}\times \int_{[0,1]^2}\int_{[0,1]^2}
\bigl(\mathfrak{f}_{q_{u}, q_{u'}}(\omega) \mathfrak {f}_{q_{v}, q_{v'}}(\omega)+ \mathfrak{f}_{q_{u}, q_{v'}}(\omega) \mathfrak{f}_{q_{v},
q_{u'}}(\omega) I\{
\omega=0 \mod\upi\} \bigr)
\\
&&\hspace*{-37pt}\hspace*{126pt}{} \times \dd \mu(u,v)\,\dd \mu\bigl(u',v'
\bigr)
\end{eqnarray*}
and the bias is given by
\[
B_{n,\mu}^{(k)}(\omega):= %
\sum_{j=2}^k \frac{b_n^j}{j!} \int
_{ -\upi}^{\upi} v^j W(v) \,\dd v
\frac
{\mathrm{d}^j}{\mathrm{d}\omega^j}\mathfrak{f}_{\mu}(\omega),\qquad \omega\in\mathbb{R}.
\]
Moreover, $Z_\mu(\omega) = Z_\mu(-\omega)$, $Z_\mu(\omega) =
Z_\mu(2\uppi+ \omega)$, and $Z_\mu(\omega)$, $\omega\in[0,\upi]$
are mutually independent random variables. The weak convergence above
holds jointly for any finite, fixed collection of frequencies $\omega$.
\end{prop}

\begin{pf}
Assumption~\textup{(C)} entails
\[
\mathfrak{f}_{\mu}(\omega) - B_{n,\mu}^{(k)}(\omega)
= \int_{[0,1]^2} \mathfrak{f}_{q_{u}, q_{v}}(\omega) -
B_n^{(k)}(u, v; \omega) \,\dd \mu(u,v).
\]
This yields
\[
\hat G_{n,\mu}(\omega) - \mathfrak{f}_{\mu}(\omega) +
B_{n,\mu
}^{(k)}(\omega) = \int_{[0,1]^2}
\mathbb{G}_{n}(u,v;\omega) \,\dd\mu(u,v),
\]
where $\mathbb{G}_{n}$ was defined in Theorem~\ref
{thmAsympDensityRankEstimator}. An application of the continuous
mapping theorem implies
\[
\sqrt{n b_n} \bigl(\hat G_{n,\mu}(\omega) -
\mathfrak{f}_{\mu
}(\omega) - B_{n,\mu}^{(k)}(\omega)
\bigr) \stackrel{\mathcal {D}} {\longrightarrow} \int_{[0,1]^2}
H(u,v;\omega) \,\dd\mu(u,v).
\]
Since $H(\cdot,\cdot;\omega)$ is a centered Gaussian process, the
integral $\int_{[0,1]^2} H(u,v;\omega) \,\dd u\,\dd v$ follows a normal
distribution with mean zero and variance:
\begin{eqnarray*}
\hspace*{-4pt}&&\int_{[0,1]^2}\int_{[0,1]^2} \Cov\bigl(H(u,v;
\omega),H\bigl(u',v';\omega \bigr)\bigr)\,\dd \mu(u,v)\,\dd\mu
\bigl(u',v'\bigr)
\\
\hspace*{-4pt}&&\quad = 2\uppi\int W^2(w) \,\dd w \int_{[0,1]^2}\int
_{[0,1]^2} \bigl( \mathfrak{f}_{q_{u}, q_{u'}}(\omega)
\mathfrak{f}_{q_{v},
q_{v'}}(\omega)+ \mathfrak{f}_{q_{u}, q_{v'}}(\omega) \mathfrak {f}_{q_{v}, q_{u'}}(
\omega) I\{\omega=0 \mod\upi\} \bigr)
\\
\hspace*{-4pt}&&\hspace*{117pt}\qquad{} \times \dd\mu(u,v)\,\dd\mu\bigl(u',v'
\bigr).
\end{eqnarray*}
This completes the proof.
\end{pf}

\subsection{The Blomqvist and Gini periodograms}
In this section, we identify two measures $\mu$ that correspond to two
classical measures of serial dependence, Blomqvist's beta (see
Blomqvist \cite{blomqvist1950},
Schmid \textit{et~al.} \cite{sstgr2010},
Genest and Carabar{\'{\i}}n-Aguirre \cite{GenestAguirreFanny2013}) and Gini's gamma
(see
Schechtman and Yitzhaki \cite{Schechtman1987},
Nelsen \cite{nelsen1998},
Carcea and Serfling \cite{Carcea2014}) coefficients, which
lead to the definition of the \textit{Blomqvist} and \textit{Gini
spectra}, respectively.

Let $C_k$ denote the copula of the pair $(X_0,X_k)$ and assume that it
is continuous. The corresponding \textit{Blomqvist~beta coefficient at lag
$k$} is
\begin{equation}
\label{beta} \beta_k:= 4 C_k(1/2,1/2) - 1.
\end{equation}
Similarly, Gini's gamma, also known as \textit{Gini's lag $k$ rank
association coefficient} is the copula-based quantity
\begin{eqnarray}\label{Gamma}
\Gamma_k &:=& 2 \int_{[0,1]^2} \bigl(\llvert u+v-1
\rrvert - \llvert v-u\rrvert \bigr) \,\dd C_k(u,v)
\nonumber\\[-8pt]\\[-8pt]\nonumber
&=& 4 \biggl( \int_{[0,1]} C_k(u,u) -
u^2 \,\dd u + \int_{[0,1]} C_k(u,1-u) -
u(1-u) \,\dd u \biggr).
\end{eqnarray}
This motivates the definition of the \textit{Blomqvist spectrum}
\[
\mathfrak{f}_{\beta}(\omega):= \frac{1}{2\uppi} \sum
_{k\in\IZ} \mathrm{e}^{-\ii\omega k} \beta_k
\]
and the \textit{Gini spectrum}
\[
\mathfrak{f}_{\Gamma}(\omega):= \frac{1}{2\uppi} \sum
_{k\in\IZ} \mathrm{e}^{-\ii\omega k} \Gamma_k.
\]
Sample versions of the Blomqvist and Gini coefficients are
\[
\hat\beta_n^k:= \frac{1}{n-\llvert  k\rrvert  }\sum
_{t=1}^{n-\llvert  k\rrvert  -1} \bigl(4 I\{ R_{n;t} \leq1/2,
R_{n;t+\llvert  k\rrvert  } \leq1/2\}- 1 \bigr),
\]
and
\[
\hat\Gamma_n^k:= \frac{2}{n(n-\llvert  k\rrvert  )} \sum
_{t=0}^{n-\llvert  k\rrvert  -1} \bigl(\llvert R_{n;t} +
R_{n;t+\llvert  k\rrvert  } - n\rrvert - \llvert R_{n;t} - R_{n;t+\llvert  k\rrvert  }
\rrvert \bigr),
\]
respectively. To establish the connection with the general periodogram
defined in the previous section, consider the measures $\mu_\beta$
which puts mass $4$ in the point $(1/2,1/2)$ and $\mu_\Gamma$ which
puts mass $4$ on the sets $\{(u,u)\dvt  u \in[0,1]\}$ and $\{(u,1-u)\dvt  u
\in[0,1]\}$, respectively.

%
%
\begin{prop}\label{propginiblomqvist}
For any $\omega\in\mathcal{F}_n$,
\[
I_{n,\beta}(\omega):= \int_{[0,1]^2} I_{n,R}^{u,v}(
\omega) \,\dd\mu _\beta(u,v) = \frac{1}{2\uppi} \sum
_{\llvert  k\rrvert   < n} \frac{n-k}{n} \ee^{\ii
\omega k} \hat
\beta_n^k
\]
and
\[
I_{n,\Gamma}(\omega):= \int_{[0,1]^2} I_{n,R}^{u,v}(
\omega) \,\dd\mu _\Gamma(u,v)= \frac{1}{2\uppi} \sum
_{\llvert  k\rrvert   < n} \frac{n-k}{n} \ee^{\ii
\omega k} \hat
\Gamma_n^k.
\]
\end{prop}

\begin{pf}
Observing that
\[
\llvert n- R_{n;t} - R_{n;t+k}\rrvert = 2 \max(n-
R_{n;t} - R_{n;t+k},0) - (n- R_{n;t} -
R_{n;t+k})
\]
and\vspace*{-4pt}
\[
\llvert R_{n;t} - R_{n;t+k}\rrvert = 2\max(R_{n;t},R_{n;t+k})
- (R_{n;t} + R_{n;t+k})
\]
yields\vspace*{-4pt}
\begin{eqnarray*}
&& \llvert R_{n;t} + R_{n;t+k} - n\rrvert -
\llvert R_{n;t} - R_{n;t+k}\rrvert
\\[-3pt]
&&\quad =2 \max(n- R_{n;t} - R_{n;t+k},0) - 2\max
(R_{n;t},R_{n;t+k}) + 2(R_{n;t} + R_{n;t+k})
- n.
\end{eqnarray*}
On the other hand,\vspace*{-3pt}
\begin{eqnarray*}
\int_{0}^1 I_{n,R}^{u,u}(
\omega)\,\mathrm{d}u &=& \frac{1}{2\uppi
}\frac
{1}{n}\sum
_{s,t=0}^{n-1} \ee^{-\ii\omega t}\ee^{\ii\omega s} \int
_{0}^1 I\{R_{n;t}\leq
nu,R_{n;s}\leq nu\} \,\mathrm{d}u
\\[-3pt]
&=& \frac{1}{2\uppi}\frac{1}{n}\sum_{s,t=0}^{n-1}
\ee^{-\ii\omega
t}\ee^{\ii\omega s} \bigl(1- n^{-1} \max(R_{n;t},R_{n;s})
\bigr)
\\[-3pt]
&=& - \frac{1}{2\uppi}\frac{1}{n^2}\sum_{s,t=0}^{n-1}
\ee^{-\ii\omega
t}\ee^{\ii\omega s} \max(R_{n;t},R_{n;s})
\end{eqnarray*}
and\vspace*{-4pt}
\begin{eqnarray*}
\int_{0}^1 I_{n,R}^{u,1-u}(
\omega)\,\mathrm{d}u &=& \frac{1}{2\uppi
}\frac{1}{n}\sum
_{s,t=0}^{n-1} \ee^{-\ii\omega t}\ee^{\ii\omega s} \int
_{0}^1 I\bigl\{R_{n;t}\leq
nu,R_{n;s}\leq n(1-u)\bigr\} \,\mathrm{d}u
\\[-3pt]
&=& \frac{1}{2\uppi}\frac{1}{n}\sum_{s,t=0}^{n-1}
\ee^{-\ii\omega
t}\ee^{\ii\omega s} \max\bigl(1 - n^{-1}R_{n;t}
- n^{-1}R_{n;s},0\bigr).
\end{eqnarray*}

Elementary algebra yields, for arbitrary functions $a$ from $\IZ^2$ to
$ \IZ$ such that $a(j,k) = a(k,j) $ for all $j,k $,\vspace*{-3pt}
\[
\sum_{\llvert  k\rrvert  <n} \sum_{t = 0}^{n-1-\llvert  k\rrvert  }
\ee^{\ii\omega k} a\bigl(t,t+\llvert k\rrvert \bigr) = \sum
_{s=0}^{n-1}\sum_{t=0}^{n-1}\ee^{-\ii\omega t}
\ee^{\ii\omega s} a(t,s).
\]
This implies (recall that $\omega\in\mathcal{F}_n$)\vspace*{-2pt}
\begin{eqnarray*}
I_{n,\Gamma}(\omega) &=& \frac{1}{2\uppi}\frac{2}{n}\sum
_{\llvert  k\rrvert  <n} \sum_{t = 0}^{n-1-\llvert  k\rrvert  }\ee^{\ii\omega k}
\bigl(\llvert R_{n;t} + R_{n;t+k} - n\rrvert - \llvert
R_{n;t} - R_{n;t+k}\rrvert \bigr)
\\[-3pt]
&=& \frac{1}{2\uppi} \frac{4}{n^2}\sum_{s=0}^{n-1}
\sum_{t=0}^{n-1}\ee^{-\ii\omega t}
\ee^{\ii\omega s} \bigl(\max(n- R_{n;t} - R_{n;s},0) -
\max(R_{n;t},R_{n;s}) \bigr)
\\[-3pt]
&&{} + \frac{1}{2\uppi} \frac{2}{n^2}\sum_{s=0}^{n-1}
\sum_{t=0}^{n-1}\ee^{-\ii\omega t}
\ee^{\ii\omega s} \bigl(2(R_{n;t} + R_{n;s}) - n \bigr)
\\[-3pt]
&=& 4 \biggl( \int_{[0,1]} I_{n,R}^{u,u}(
\omega) \,\dd u + \int_{[0,1]} I_{n,R}^{u,1-u}(
\omega) \,\dd u \biggr).
\end{eqnarray*}
The representation for $I_{n,\beta}$ can be derived similarly; details
are omitted for the sake of brevity.
\end{pf}

Smoothed versions of the Blomqvist and Gini periodograms can be defined
accordingly, and their asymptotic distributions follow from
Proposition~\ref{propAsympGenLin}. In particular, this yields
consistent estimators of the Blomqvist and Gini spectra defined above.

We conclude this section with some general remarks. First, note that
the approach above can be applied to any scalar dependence measure that
can be represented as a continuous linear functional of the copula. For
instance, Cifarelli, Conti and Regazzini \cite{CifarelliContiRegazzini1996} consider a general
measure of monotone dependence of the form
\begin{equation}
\label{eqgencon} \int_{[0,1]^2} g\bigl(\llvert u+v-1\rrvert \bigr)
- g\bigl(\llvert u-v\rrvert \bigr) \,\dd C(u,v),
\end{equation}
where $g\dvtx  [0,1] \to\mathbb{R}$ is strictly increasing and convex.
Choosing $g(x) = x$ and $g(x) = x^2$ yields (up to constants) the Gini
and Spearman rank correlations, respectively. Under suitable
assumptions on $g$, the monotone dependence measure in (\ref
{eqgencon}) can be written (by applying integration-by-parts) in the
form of equation (\ref{eqinmu}), and the results from Section~\ref{secgenlin} apply.

Other measures of serial dependence such as Kendall's $\tau$ (see
Ferguson, Genest and Hallin \cite{FergusonEtAl2000}) only can
be represented as nonlinear
functionals of the copula. More general rank-based autocorrelation
coefficients also have been introduced in the context of inference for
ARMA models (see Hallin and Puri \cite{HallinPuri1992,Hallin1994} or Hallin
\cite{Hallin2012}); they involve score functions, typically are not
time-revertible, and lead to possibly unbounded measures $\mu$. We
expect that the general results presented here can be extended to the
periodograms associated with such coefficients, but leave this question
to future research.

\section{Simulation study} \label{secsim}

In this section, we show how Theorem~\ref
{thmAsympDensityRankEstimator} can be used to construct asymptotic
confidence intervals for the copula spectra. An analysis of the finite
sample performance was conducted using the \proglang{R} package \pkg
{quantspec} (Kley \cite{Kley2014,Kley2014a}). We consider three different models:
\begin{longlist}[(a)]
\item[(a)] the QAR(1)~process
\begin{equation}
\label{eqnmodel1} Y_t = 0.1 \Phi^{-1}(U_t) +
1.9(U_t-0.5) Y_{t-1}
\end{equation}
(cf.~Koenker and Xiao \cite{Koenker2006}), where $(U_t)$
is a sequence of i.i.d.
standard uniform random variables, and~$\Phi$ denotes the distribution
function of the standard normal distribution;
\item[(b)] the AR(2)~process
\begin{equation}
\label{eqnmodel2} Y_t = - 0.36 Y_{t-2} +
\varepsilon_t,
\end{equation}
where $(\varepsilon_t)$ is standard normal white noise
(cf.~Li \cite{Li2012});
\item[(c)] the ARCH(1)~process
\begin{equation}
\label{eqnmodel3} Y_t = \bigl( 1/1.9 + 0.9 Y_{t-1}^2
\bigr)^{1/2} \varepsilon_t,
\end{equation}
where $(\varepsilon_t)$ is standard normal white noise
(cf.~Lee and Rao \cite{LeeRao2012}).
\end{longlist}

For each model, 10\,000 independent copies of length $n \in\{2^8,
2^9, 2^{10}, 2^{11}\}$ were generated.
For each of them, the smoothed CR-periodograms
\begin{equation}
\label{normalized} \tilde G_{n,R}(\tau_1, \tau_2;
\omega_{jn}):= {\hat G_{n,R}(\tau _1,
\tau_2; \omega_{jn})}/{W_n^j},
\qquad W_n^j:= \frac{2\uppi}{n} \sum
_{{0=s \neq j}}^{n-1} W_n ( \omega _{jn}
- \omega_{sn} ),
\end{equation}
were computed for $\omega_{jn}:= 2\uppi j/n$, $j=1,\ldots,n/2-1$ and
$\tau_1, \tau_2 \in\{0.1, 0.5, 0.9\}$, where we used the kernel of
order $4$
\[
W(u):= \frac{15}{32} \frac{1}{\upi} \bigl(7 ({u}/{\upi})^4
- 10 ({u}/{\upi})^2 + 3 \bigr) I\bigl\{\llvert u \rrvert \leq\upi\bigr
\}
\]
minimizing the asymptotic IMSE (see
Gasser, M{\"u}ller and Mammitzsch
\cite{GasserEtAl1985}). The
bandwidth was chosen as $b_n = 0.4 n^{-1/4}$ which is of lower order
than the IMSE-optimal bandwidth~$ n^{-1/9} $ to reduce bias and the
factor $(W_n^j)^{-1}$ ensures that the weights in (\ref{normalized})
sum up to one for every~$n$.

Based on Theorem~\ref{thmAsympDensityEstimator}, we then computed
pointwise asymptotic $(1-\alpha)$-level confidence bands for the real
and imaginary parts of the spectrum, namely,
\begin{equation}
\label{eqnDefCI} IC_{1,n}(\tau_1, \tau_2;
\omega_{jn}):=\Re\tilde G_{n,R}(\tau_1, \tau _2;
\omega_{jn}) \pm\Re\sigma(\tau_1, \tau_2; \omega_{jn}) \Phi
^{-1}(1-\alpha/2),
\end{equation}
for the real part, and
\begin{equation}
\label{eqnDefCI-im} IC_{2,n}(\tau_1, \tau_2;
\omega_{jn}):= \Im\tilde G_{n,R}(\tau_1, \tau _2;
\omega_{jn}) \pm\Im\sigma(\tau_1, \tau_2; \omega_{jn}) \Phi
^{-1}(1-\alpha/2),
\end{equation}
for the imaginary part of the copula spectrum. As usual, $\Phi$ stands
for the standard normal distribution function, and
\[
\bigl( \Re\sigma(\tau_1, \tau_2; \omega_{jn})
\bigr)^2:= 0 \vee %
\cases{ \displaystyle c(
\tau_1, \tau_2; \omega_{jn}, \omega_{jn}), &\quad if $
\tau_1 = \tau_2$,
\vspace*{3pt}\cr
\displaystyle\tfrac{1}{2}
\bigl(c(\tau_1, \tau_2; \omega_{jn}, \omega_{jn}) + c(
\tau_1, \tau_2; \omega_{jn}, -\omega_{jn}) \bigr), &\quad if $
\tau_1 \neq\tau_2$,}
\]
and
\[
\bigl( \Im\sigma(\tau_1, \tau_2; \omega_{jn})
\bigr)^2:= 0 \vee %
\cases{ 0, &\quad if $\tau_1
= \tau_2$,
\vspace*{3pt}\cr
\displaystyle\tfrac{1}{2} \bigl(c(
\tau_1, \tau_2; \omega_{jn}, \omega_{jn}) - c(\tau_1,
\tau_2; \omega_{jn}, -\omega_{jn}) \bigr), &\quad if $\tau_1 \neq
\tau_2$}
\]
are estimators for $\Var (\Re\tilde G_{n,R}(\tau_1, \tau_2;
\omega_{jn})  )$ and $\Var (\Im\tilde G_{n,R}(\tau_1, \tau_2;
\omega_{jn})  )$, respectively. Here,
\begin{eqnarray*}
&& c\bigl(\tau_1, \tau_2; \omega_{jn}, \omega_{j'n}
\bigr)
\\
&&\quad := \biggl(\frac{2\uppi}{n} \Big/W_n^j
\biggr)^2
\\
&&\qquad{}\times \Biggl[ \sum_{s=1}^{n-1}
W_n (\omega_{jn}- 2\uppi s / n ) W_n \bigl(\omega_{j'n}
- 2\uppi s / n \bigr) \tilde G_{n,R}(\tau_1,
\tau_1; 2\uppi s / n) \tilde G_{n,R}(\tau_2,
\tau_2; 2\uppi s / n)
\\
&&\hspace*{98pt}{}+ \sum_{s=1}^{n-1} W_n (\omega_{jn}-
2\uppi s / n ) W_n \bigl(\omega_{j'n} + 2\uppi s / n \bigr)
\bigl\llvert \tilde G_{n,R}(\tau_1, \tau_2; 2\uppi s
/ n) \bigr\rrvert ^2 \Biggr]
\end{eqnarray*}
is an estimator for the covariance of $\tilde G_{n,R}(\tau_1, \tau_2;
\omega_{jn})$ and $\tilde G_{n,R}(\tau_1, \tau_2; \omega_{j'n})$; this
follows from the representation in Theorem~\ref
{thmAsympDensityEstimator}(iii) and Theorem~7.4.3 in Brillinger \cite
{Brillinger1975}. To motivate this approach,
recall that, for any
complex-valued random variable $Z$ with complex conjugate $\bar{Z}$,
\[
\Var(\Re Z) = \tfrac{1}{2} \bigl( \Var(Z) + \Re\Cov(Z, \bar {Z}) \bigr);\qquad
\Var(\Im Z) = \tfrac{1}{2} \bigl( \Var(Z) - \Re\Cov(Z, \bar {Z}) \bigr).
\]
For $n \to\infty$, the estimated variances above converge to the
asymptotic variance in Theorem \ref{thmAsympDensityRankEstimator}.
However, in small samples the more elaborate version considered here
typically leads to better coverage probabilities.

%
%
\begin{table}[b]
\tabcolsep=0pt
\caption{Coverage frequencies for the confidence intervals $IC_n(\tau
_1, \tau_2, \omega)$, $n = 2^8$, $b_n = 0.4 n^{-1/4}$, $1-\alpha= 0.95$}\label{tabCI-2h8}
\begin{tabular*}{\tablewidth}{@{\extracolsep{\fill}}@{}lllllll@{}}
\hline
& $(\tau_1, \tau_2)$ & $(0.1,0.1)$ & $(0.1,0.9)$ & $(0.5,0.5)$ & $(0.1,0.9)$ & $(0.9,0.9)$ \\
Model & $\omega/\upi$ & $(\Re)$ & $(\Im)$ & $(\Re)$ & $(\Re)$ & $(\Re)$ \\
\hline
(a) QAR(1) (\ref{eqnmodel1}) & $1/8$ & 0.911 & 0.921 & 0.906 & 0.987
& 0.899 \\
& $1/4$ & 0.934 & 0.917 & 0.920 & 0.979 & 0.910 \\
& $1/2$ & 0.947 & 0.919 & 0.932 & 0.976 & 0.915 \\
& $3/4$ & 0.946 & 0.918 & 0.927 & 0.979 & 0.916 \\
& $7/8$ & 0.941 & 0.915 & 0.931 & 0.979 & 0.921
\\[3pt]
(b) AR(2) (\ref{eqnmodel2}) & $1/8$ & 0.913 & 0.926 & 0.900 & 0.975 &
0.916 \\
& $1/4$ & 0.935 & 0.925 & 0.917 & 0.967 & 0.940 \\
& $1/2$ & 0.940 & 0.927 & 0.929 & 0.966 & 0.949 \\
& $3/4$ & 0.939 & 0.924 & 0.928 & 0.969 & 0.947 \\
& $7/8$ & 0.937 & 0.920 & 0.928 & 0.972 & 0.945
\\[3pt]
(c) ARCH(1) (\ref{eqnmodel3}) & $1/8$ & 0.860 & 0.910 & 0.906 & 0.902
& 0.878 \\
& $1/4$ & 0.872 & 0.905 & 0.922 & 0.909 & 0.887 \\
& $1/2$ & 0.902 & 0.897 & 0.937 & 0.946 & 0.914 \\
& $3/4$ & 0.906 & 0.894 & 0.934 & 0.959 & 0.924 \\
& $7/8$ & 0.906 & 0.891 & 0.935 & 0.962 & 0.920 \\
\hline
\end{tabular*}
\end{table}

%
%
\begin{table}
\tabcolsep=0pt
\caption{Coverage frequencies for the confidence intervals $IC_n(\tau
_1, \tau_2, \omega)$, $n = 2^9$, $b_n = 0.4 n^{-1/4}$,
$1-\alpha= 0.95$}\label{tabCI-2h9}%
\begin{tabular*}{\tablewidth}{@{\extracolsep{\fill}}@{}lllllll@{}}
\hline
& $(\tau_1, \tau_2)$ & $(0.1,0.1)$ & $(0.1,0.9)$ & $(0.5,0.5)$ & $(0.1,0.9)$ & $(0.9,0.9)$ \\
Model & $\omega/\upi$ & $(\Re)$ & $(\Im)$ & $(\Re)$ & $(\Re)$ & $(\Re)$ \\
\hline
(a) QAR(1) (\ref{eqnmodel1}) & $1/8$ & 0.934 & 0.932 & 0.915 & 0.974
& 0.916 \\
& $1/4$ & 0.953 & 0.933 & 0.931 & 0.968 & 0.925 \\
& $1/2$ & 0.954 & 0.932 & 0.940 & 0.968 & 0.934 \\
& $3/4$ & 0.952 & 0.926 & 0.939 & 0.973 & 0.932 \\
& $7/8$ & 0.953 & 0.923 & 0.941 & 0.975 & 0.934
\\[3pt]
(b) AR(2) (\ref{eqnmodel2}) & $1/8$ & 0.930 & 0.934 & 0.913 & 0.962 &
0.932 \\
& $1/4$ & 0.950 & 0.932 & 0.928 & 0.956 & 0.951 \\
& $1/2$ & 0.948 & 0.935 & 0.933 & 0.957 & 0.949 \\
& $3/4$ & 0.951 & 0.932 & 0.936 & 0.964 & 0.952 \\
& $7/8$ & 0.949 & 0.931 & 0.937 & 0.965 & 0.955
\\[3pt]
(c) ARCH(1) (\ref{eqnmodel3}) & $1/8$ & 0.890 & 0.932 & 0.918 & 0.913 & 0.892 \\
& $1/4$ & 0.900 & 0.924 & 0.938 & 0.917 & 0.903 \\
& $1/2$ & 0.922 & 0.912 & 0.939 & 0.948 & 0.928 \\
& $3/4$ & 0.926 & 0.913 & 0.944 & 0.957 & 0.934 \\
& $7/8$ & 0.928 & 0.908 & 0.943 & 0.958 & 0.937 \\
\hline
\end{tabular*}
\end{table}
%

%
%
\begin{table}
\tabcolsep=0pt
\caption{Coverage frequencies for the confidence intervals $IC_n(\tau
_1, \tau_2, \omega)$, $n = 2^{10}$, $b_n = 0.4 n^{-1/4}$,
$1-\alpha= 0.95$}\label{tabCI-2h10}%
\begin{tabular*}{\tablewidth}{@{\extracolsep{\fill}}@{}lllllll@{}}
\hline
& $(\tau_1, \tau_2)$ & $(0.1,0.1)$ & $(0.1,0.9)$ & $(0.5,0.5)$ & $(0.1,0.9)$ & $(0.9,0.9)$ \\
Model & $\omega/\upi$ & $(\Re)$ & $(\Im)$ & $(\Re)$ & $(\Re)$ & $(\Re)$ \\
\hline
(a) QAR(1) (\ref{eqnmodel1}) & $1/8$ & 0.942 & 0.943 & 0.933 & 0.961
& 0.924 \\
& $1/4$ & 0.959 & 0.938 & 0.941 & 0.963 & 0.929 \\
& $1/2$ & 0.953 & 0.938 & 0.941 & 0.962 & 0.934 \\
& $3/4$ & 0.954 & 0.935 & 0.941 & 0.967 & 0.933 \\
& $7/8$ & 0.956 & 0.935 & 0.943 & 0.969 & 0.936
\\[3pt]
(b) AR(2) (\ref{eqnmodel2}) & $1/8$ & 0.939 & 0.943 & 0.931 & 0.953 &
0.940 \\
& $1/4$ & 0.954 & 0.939 & 0.942 & 0.954 & 0.952 \\
& $1/2$ & 0.954 & 0.944 & 0.945 & 0.953 & 0.955 \\
& $3/4$ & 0.950 & 0.937 & 0.942 & 0.956 & 0.954 \\
& $7/8$ & 0.954 & 0.937 & 0.940 & 0.959 & 0.952
\\[3pt]
(c) ARCH(1) (\ref{eqnmodel3}) & $1/8$ & 0.900 & 0.935 & 0.933 & 0.911 & 0.906 \\
& $1/4$ & 0.901 & 0.930 & 0.945 & 0.916 & 0.908 \\
& $1/2$ & 0.929 & 0.928 & 0.945 & 0.942 & 0.928 \\
& $3/4$ & 0.941 & 0.916 & 0.948 & 0.954 & 0.937 \\
& $7/8$ & 0.940 & 0.918 & 0.948 & 0.953 & 0.936 \\
\hline
\end{tabular*}
\end{table}
%

%
%
\begin{table}
\tabcolsep=0pt
\caption{Coverage frequencies for the confidence intervals $IC_n(\tau
_1, \tau_2, \omega)$, $n = 2^{11}$, $b_n = 0.4 n^{-1/4}$, $1-\alpha= 0.95$}\label{tabCI-2h11}%
\begin{tabular*}{\tablewidth}{@{\extracolsep{\fill}}@{}@{}lllllll@{}}
\hline
& $(\tau_1, \tau_2)$ & $(0.1,0.1)$ & $(0.1,0.9)$ & $(0.5,0.5)$ & $(0.1,0.9)$ & $(0.9,0.9)$ \\
Model & $\omega/\upi$ & $(\Re)$ & $(\Im)$ & $(\Re)$ & $(\Re)$ & $(\Re)$ \\
\hline
(a) QAR(1) (\ref{eqnmodel1}) & 1/8 & 0.953 & 0.945 & 0.944 & 0.957 & 0.933 \\
& 1/4 & 0.957 & 0.943 & 0.945 & 0.961 & 0.932 \\
& 1/2 & 0.955 & 0.938 & 0.949 & 0.960 & 0.938 \\
& 3/4 & 0.952 & 0.938 & 0.946 & 0.963 & 0.939 \\
& 7/8 & 0.954 & 0.936 & 0.945 & 0.964 & 0.945
\\[3pt]
(b) AR(2) (\ref{eqnmodel2}) & 1/8 & 0.953 & 0.944 & 0.943 & 0.954 &
0.947 \\
& 1/4 & 0.954 & 0.944 & 0.945 & 0.953 & 0.956 \\
& 1/2 & 0.955 & 0.946 & 0.945 & 0.951 & 0.954 \\
& 3/4 & 0.954 & 0.947 & 0.940 & 0.954 & 0.957 \\
& 7/8 & 0.952 & 0.945 & 0.943 & 0.956 & 0.951
\\[3pt]
(c) ARCH(1) (\ref{eqnmodel3}) & 1/8 & 0.911 & 0.942 & 0.944 & 0.918 & 0.908 \\
& 1/4 & 0.918 & 0.937 & 0.950 & 0.926 & 0.917 \\
& 1/2 & 0.934 & 0.931 & 0.947 & 0.946 & 0.937 \\
& 3/4 & 0.944 & 0.931 & 0.949 & 0.954 & 0.943 \\
& 7/8 & 0.944 & 0.928 & 0.950 & 0.958 & 0.945 \\
\hline
\end{tabular*}
\end{table}

In Tables~\ref{tabCI-2h8}--\ref{tabCI-2h11}, we report the
simulated coverage frequencies associated with
\[
\IP \bigl(\Re\mathfrak{f}_{q_{\tau_1}, q_{\tau_2}}(\omega) \in IC_{1,n} (
\tau_1, \tau_2, \omega) \bigr)\quad\mbox{and}\quad \IP
\bigl(\Im\mathfrak{f}_{q_{\tau_1}, q_{\tau_2}}(\omega) \in IC_{2,n} (
\tau_1, \tau_2, \omega) \bigr).
\]

Inspection of Tables~\ref{tabCI-2h8}--\ref{tabCI-2h11} reveals
that, as $n$ gets larger, the coverage frequencies converge to the
confidence level $1-\alpha$. For models (\ref{eqnmodel1})--(\ref
{eqnmodel2}), those frequencies are quite close to $1-\alpha$ even for
moderately large values of~$n$. Due to boundary effects, the coverage
frequencies for $\omega$ close to multiples of $\upi$ are too low in
all three models, but, as noted earlier, they improve as $n$ increases.
Finally, in models~(\ref{eqnmodel1}) and~(\ref{eqnmodel3}) for
smaller values of~$n$, the confidence intervals involving extreme
quantiles tend to cover less frequently, as can be expected. Again, the
accuracy improves with increasing sample size.

\section{Conclusions}\label{secConcl}

Spectral analysis for the past fifty years has been a major tool in the
analysis of time series. Being essentially covariance-based, however,
classical $L^2$-spectral methods have obvious limitations, for instance
(see Figures~\ref{figL2Spectra} and~\ref{figCopulaSpectra}), they
cannot discriminate between QAR or ARCH and white noise processes.
Quantile-related spectral concepts have been proposed, which palliate
those limitations. Only quite incomplete asymptotic distributional
results, however, have been available in the literature for the
consistent estimation of such concepts, which so far has precluded most
practical applications.

In this paper, we provide (Theorem~\ref
{thmAsympDensityRankEstimator}), in the very strong form of
convergence to a Gaussian process, such asymptotic results for the
smoothed copula rank-based periodogram process. That rank-based
periodogram is the generalization (Dette \textit{et~al.}
\cite{DetteEtAl2013}) of the copula
rank periodograms proposed by Hagemann \cite
{Hagemann2013}. Theorem~\ref
{thmAsympDensityRankEstimator} was used to construct confidence
intervals. A simulation study was conducted using the \proglang{R}
package \pkg{quantspec} (Kley \cite{Kley2014,Kley2014a}).

Being copula- or rank-based, our spectral concepts furthermore are
invariant under monotone increasing continuous marginal transformations
of the data, and are likely to enjoy appealing robustness features
their traditional $L^2$-counterparts are severely lacking. Another
application is in the asymptotic behavior of the spectra associated
with more classical rank-based autocorrelation coefficients, such as
the Spearman, Gini or Blomqvist spectra.

Copula rank-based periodogram methods are improving over the classical
ones both from the point of view of efficiency (detection of nonlinear
features) and from the point of view of robustness (no finite variance
assumption is required). They are likely to be ideal tools for a large
variety of problems of practical interest, such as change-point
analysis, tail dependence, model diagnostics, or local stationary
procedures~(cf. Skowronek \cite{Skowronek2014})~-- essentially, all problems
covered in the traditional spectral context can be extended here, with
the huge advantage that nonlinear features that cannot be accounted for
by traditional methods can be analyzed via the new ones. This seems to
offer most promising perspectives for future research.

\begin{appendix}\label{pmt}
\section*{Appendix: Proof of Theorem \texorpdfstring{\protect\ref{thmAsympDensityEstimator}}{3.6}}

The proof of Theorem~\ref{thmAsympDensityEstimator} relies on a
series of technical lemmas; for the readers' convenience, we begin by
giving a general overview of the main steps and the corresponding lemmas.

For all $n\in\IN$, consider the stochastic process
\begin{equation}
\label{eqdefhn} \hat H_{n,U}(\tau_1,\tau_2;
\omega):= \sqrt{nb_n} \bigl(\hat G_{n,U}(
\tau_1,\tau_2;\omega) - \IE\hat G_{n,U}(
\tau_1,\tau _2;\omega) \bigr),
\end{equation}
indexed by $(\tau_1,\tau_2)\in[0,1]^2$ and $\omega\in\IR$; for $a
= (a_1,a_2)\in[0,1]^2$, write~$\hat H_n(a;\omega)$ for~$\hat
H_{n,U}(a_1, a_2;\omega)$.

The key step in the process of establishing parts (i) and (iii) of
Theorem~\ref{thmAsympDensityEstimator} is a uniform bound on the
increments of the process $\hat H_{n,U}$. That bound is required, for
example, when showing the stochastic equicontinuity of $\hat
H_n(a;\omega) - \hat H_n(b;\omega)$. We derive such a bound by a
restricted chaining technique, which is described in Lemma~\ref
{lemThm224VW}. The application of Lemma~\ref{lemThm224VW} requires
two ingredients. First, we need a general bound, uniform in $a$ and
$b$, on the moments of $\hat H_n(a;\omega) - \hat H_n(b;\omega)$.
Such a bound is derived in Lemma~\ref{lemSixthMomIncrementHn}.
Second, we need a sharper bound on the increments $\hat H_n(a;\omega)
- \hat H_n(b;\omega)$ when $a$ and $b$ are ``sufficiently close''. We
provide this result in Lemma~\ref{lemOrderSmallIncrements}.

Lemma~\ref{lemSixthMomIncrementHn} is a very general result, relying
on an abstract condition on the cumulants of discrete Fourier
transforms of certain indicator functions; see~(\ref
{lemSixthMomIncrementHneqnAssumption}). The link between assumption~\textup{(C)} and~(\ref{lemSixthMomIncrementHneqnAssumption}) is established
in Lemma~\ref{lemOrderDftYinA}.

Finally, the proof of part (ii) of Theorem~\ref
{thmAsympDensityEstimator} follows by a series of uniform
generalizations of results from~Brillinger \cite{Brillinger1975}, the details of
which are provided in the online supplement~\cite{supp} [Lemmas 8.1--8.5].

\subsection{Proof of part~\textup{(i)} of Theorem \texorpdfstring{\protect\ref{thmAsympDensityEstimator}}{3.6}}\label{Parti}
In view of Theorems~1.5.4
and~1.5.7 in van~der Vaart and Wellner \cite{vanderVaartWellner1996}, it is sufficient to prove
the following two claims:
\begin{longlist}[(i2)]
\item[(i1)] convergence of the finite-dimensional distributions of the
process (\ref{eqdefhn}), that is,
\begin{equation}
\label{thmAsympDensityEstimatoreqnfidis} \bigl(\hat H_n(a_{1j},a_{2j};
\omega_j) \bigr)_{j=1,\ldots,k} \xrightarrow{d} \bigl(H(a_{1j},a_{2j};
\omega_j) \bigr)_{j=1,\ldots,k}
\end{equation}
for any $(a_{1j},a_{2j},\omega_j) \in[0,1]^2\times\IR$, $j=1,\ldots,k$ and $k \in\IN$;

\item[(i2)] stochastic equicontinuity:
for any $x > 0$ and any $\omega\in\IR$,
\begin{equation}
\label{thmAsympDensityEstimatoreqnstochequicont} \lim_{\delta\downarrow0} \limsup_{n \rightarrow\infty}
\IP \Bigl( \mathop{\sup_{a, b \in[0,1]^2}}_{\llVert   a - b \rrVert  _1 \leq\delta} \bigl\llvert
\hat H_n(a;\omega) - \hat H_n(b; \omega) \bigr\rrvert > x
\Bigr) = 0.
\end{equation}
\end{longlist}
Note indeed that~(\ref{thmAsympDensityEstimatoreqnstochequicont})
implies stochastic equicontinuity of both the real part~$ ( \Re
\hat H_n(a;\break \omega)  )_{a \in[0,1]^2}$ and the imaginary part
$ ( \Im\hat H_n(a;\omega)  )_{a \in[0,1]^2}$ of~$\hat H_n$.\vspace*{2pt}

First consider (i1).
Observe that $\hat G_{n,U}(\tau_1,\tau_2;\omega)$ is the traditional
smoothed periodogram estimator (see Chapter~7.1 in~Brillinger \cite{Brillinger1975}) of the cross-spectrum of
the \textit{clipped processes}
$(I\{F(X_t)\leq\tau_1\})_{t\in\IZ}$ and $ (I\{F(X_t)\leq\tau_2\}
)_{t\in\IZ}$. Thus,~(\ref{thmAsympDensityEstimatoreqnfidis}) is
an immediate corollary of Theorem~7.4.4 in~Brillinger
\cite{Brillinger1975}. The
limiting first and second moment structures are given by Theorem~7.4.1
and Corollary~7.4.3 in~Brillinger \cite
{Brillinger1975}. This implies the desired
convergence (\ref{thmAsympDensityEstimatoreqnfidis}) of
finite-dimensional distributions. Note that, by condition~\textup{(C)}, the
summability condition required for the three theorems holds
(Assumption~2.6.2($\ell$), for every $\ell$; cf.~Brillinger \cite
{Brillinger1975}).

Turning to (i2), in the notation from van~der Vaart and Wellner  \cite{vanderVaartWellner1996},
page 95, put $\Psi(x):= x^6$: the Orlicz norm $\llVert   X \rrVert  _{\Psi} =
\inf\{C > 0\dvt  \IE\Psi( \llvert  X \rrvert   / C ) \leq1\}$ coincides with the $L_6$
norm~$\llVert   X \rrVert  _6 = ( \IE\llvert  X \rrvert  ^6 )^{1/6}$. Therefore, by Lemma~\ref
{lemSixthMomIncrementHn} and Lemma~\ref{lemOrderDftYinA}, we have,
for any $\kappa\in (0,1)$ and sufficiently small $\llVert  a-b\rrVert  _1$,
\[
\bigl\llVert \hat H_n(a;\omega) - \hat H_n(b;\omega)
\bigr\rrVert _{\Psi} \leq K \biggl( \frac{\llVert  a-b\rrVert  _1^\kappa}{(n b_n)^2} +
\frac{\llVert  a-b\rrVert  _1^{2\kappa}}{n
b_n} + \llVert a-b\rrVert _1^{3\kappa}
\biggr)^{1/6}.
\]
It follows that, for all $a, b$ with $\llVert  a-b\rrVert  _1$ sufficiently small and
$\llVert  a-b\rrVert  _1 \geq(n b_n)^{-1/\gamma}$ and all $\gamma\in(0,1)$ such
that $\gamma< \kappa$,
\begin{eqnarray*}
\bigl\llVert \hat H_n(a;\omega) - \hat H_n(b;\omega)
\bigr\rrVert _{\Psi} & \leq&K \bigl( \llVert a-b\rrVert
_1^{\kappa+2\gamma} + \llVert a-b\rrVert _1^{2\kappa
+\gamma}
+ \llVert a-b\rrVert _1^{3\kappa} \bigr)^{1/6}
\\
& \leq& \bar K \llVert a-b\rrVert _1^{\gamma/2}.
\end{eqnarray*}
Note that\vspace*{1pt}
$
\llVert  a-b\rrVert  _1\geq(n b_n)^{-1/\gamma}$ iff
$d(a,b):= \llVert  a-b\rrVert  _1^{\gamma/2} \geq(n b_n)^{-1/2} =: \bar\eta_n /
2$.

Denoting by $D(\varepsilon,d)$ the \textit{packing number} of
$([0,1]^2,d)$ (cf. van~der Vaart and Wellner \cite{vanderVaartWellner1996}, page~98), we
have~$D(\varepsilon,d) \asymp\varepsilon^{-4/\gamma}$. Therefore,
by Lemma~\ref{lemThm224VW}, for all $x, \delta> 0$ and $\eta\geq
\bar\eta_n$,
\begin{eqnarray*}
&& \IP \Bigl( \sup_{\llVert   a - b\rrVert  _1 \leq\delta^{2/\gamma}} \bigl\llvert \hat H_n(a;
\omega) - \hat H_n(b;\omega)\bigr\rrvert > x \Bigr)
\\
&&\quad = \IP \Bigl( \sup_{d(a,b) \leq\delta} \bigl\llvert \hat
H_n(a;\omega) - \hat H_n(b;\omega)\bigr\rrvert > x \Bigr)
\\
&&\quad \leq \biggl[ \frac{8 \tilde K}{x} \biggl(\int_{\bar\eta
_n/2}^{\eta}
\epsilon^{-2/(3\gamma)} \,\mathrm{d}\epsilon+ (\delta+ 2 \bar\eta_n)
\eta^{-4/(3\gamma)} \biggr) \biggr]^6
\\
&&\qquad{} + \IP \Bigl( \sup_{d(a,b) \leq\bar\eta
_n} \bigl\llvert
\hat H_n(a;\omega) - \hat H_n(b;\omega) \bigr\rrvert >
x/4 \Bigr).
\end{eqnarray*}

Now choose $1>\gamma>2/3$. Letting $n$ tend to infinity, the second
term tends to zero by Lemma~\ref{lemOrderSmallIncrements} since, by
construction, $1/\gamma> 1$ and
\[
d(a,b) \leq\bar\eta_n\quad\mbox{iff}\quad \llVert a-b\rrVert
_1 \leq2^{2/\gamma} (n b_n)^{-1/\gamma}.
\]

All together, this implies
\[
\lim_{\delta\downarrow0} \limsup_{n \rightarrow\infty} \IP \Bigl( \sup
_{d(a,b) \leq\delta} \bigl\llvert \hat H_n(a;\omega) - \hat
H_n(b;\omega )\bigr\rrvert > x \Bigr) \leq \biggl[ \frac{8 \tilde K}{x}
\int_{0}^{\eta} \epsilon ^{-2/(3\gamma)} \,\mathrm{d}
\epsilon \biggr]^6,
\]
for every $x, \eta> 0$; the claim follows, since the integral in the
right-hand side can be made arbitrarily small by choosing $\eta$
accordingly.

\subsection{Proof of part~\textup{(ii)} of Theorem \texorpdfstring{\protect\ref{thmAsympDensityEstimator}}{3.6}}\label{Partii} Essentially, this part of
Theorem~\ref{thmAsympDensityEstimator} is a uniform version of
Theorems~7.4.1 and 7.4.2 in~Brillinger \cite
{Brillinger1975} in the present
setting of Laplace spectra. The proof is based on a series of uniform
versions of results from~Brillinger \cite
{Brillinger1975}; details are provided
in the online supplement \cite{supp}
(see in particular Lemma~8.5).

\subsection{Proof of part~\textup{(iii)} of Theorem \texorpdfstring{\protect\ref{thmAsympDensityEstimator}}{3.6}}\label{Partiii}
It follows from the continuity of~$F$ that the ranks of the random
variables $X_0,\ldots,X_{n-1}$ and $F(X_0),\ldots,F(X_{n-1})$
coincide almost
surely. Thus, without loss of generality, we can assume that the
estimator is computed from the unobservable data
$F(X_0),\ldots,F(X_{n-1})$. In particular, this implies that we can assume
the marginals to be uniform.

Denote by $\hat F_n^{-1}(\tau):= \inf\{x\dvt  \hat{F}_n(x) \geq\tau\}$ the
generalized inverse of $\hat F_n$ and let $\inf\varnothing:= 0$.
Elementary computation shows that, for any $k \in\IN$,
\begin{equation}
\label{eqdntau} \sup_{\omega\in\IR}\sup_{\tau\in[0,1]}\bigl
\llvert d_{n,R}^{\tau
}(\omega) - d_{n,U}^{\hat F_n^{-1}(\tau)}(
\omega)\bigr\rrvert \leq n \sup_{\tau\in[0,1]} \bigl\llvert \hat
F_n(\tau) - \hat F_n(\tau-)\bigr\rrvert =
\mathrm{O}_P\bigl(n^{1/2k}\bigr),
\end{equation}
where $\hat F_n(\tau-):= \lim_{\xi\uparrow0} \hat F_n(\tau-\xi)$
and the $\mathrm{O}_P$-bound in the above equation follows from Lemma~8.6 (online supplement~\cite{supp}). By the definition of $\hat G_{n,R}$ and arguments
similar to the ones used in the proof of Lemma~\ref
{lemOrderSmallIncrements}, it follows that
\[
\sup_{\omega\in\IR}\sup_{\tau_1, \tau_2 \in[0,1]} \bigl\llvert \hat
G_{n,R} (\tau_1, \tau_2; \omega) - \hat
G_{n,U} \bigl(\hat F_n^{-1}(\tau _1),
\hat F_n^{-1}(\tau_2); \omega\bigr) \bigr
\rrvert = \mathrm{o}_P(1).
\]
It therefore suffices to bound the differences
\[
\sup_{\tau_1, \tau_2 \in[0,1]} \bigl\llvert \hat G_{n,U} (
\tau_1, \tau_2; \omega) - \hat G_{n,U} \bigl(
\hat F_n^{-1}(\tau_1), \hat F_n^{-1}(
\tau _2); \omega\bigr) \bigr\rrvert
\]
pointwise and uniformly in $\omega$.

In what follows, we give a detailed proof of the statement for fixed
$\omega\in\IR$ and sketch the arguments needed for the proof of the
uniform result.

By~(\ref{eqdefhn}) we have, for any $x>0$ and $\delta_n$ with
\[
n^{-1/2} \ll\delta_n = \mathrm{o}\bigl( n^{-1/2}
b_n^{-1/2} (\log n)^{-d} \bigr),
\]
where $d$ is the constant from Lemma~\ref{lemLipschitzLaplaceSD}
corresponding to $j=k$,
\begin{eqnarray*}
&& P^n(\omega)
\\
&&\quad := \IP \Bigl( \sup_{\tau_1, \tau_2 \in[0,1]} \bigl
\llvert \hat G_{n,U} \bigl(\hat F_n^{-1}(
\tau_1), \hat F_n^{-1}(\tau_2);
\omega\bigr) - \hat G_{n,U} (\tau_1, \tau_2;
\omega) \bigr\rrvert > x \bigl((nb_n)^{-1/2}+b_n^k
\bigr) \Bigr)
\nonumber
\\
&&\quad \leq \IP \Bigl( \sup_{\tau_1, \tau_2 \in[0,1]} \mathop{\sup_{\llVert
(u,v) - (\tau_1, \tau_2)\rrVert  _{\infty}}}_{{\leq\sup_{\tau\in
[0,1]} \llvert   \hat F_n^{-1}(\tau) - \tau\rrvert   }}
\bigl\llvert \hat G_{n, U} ( u, v; \omega) - \hat G_{n,U} (
\tau_1, \tau_2; \omega) \bigr\rrvert > x
\bigl((nb_n)^{-1/2}+b_n^k\bigr) \Bigr)
\nonumber
\\
&&\quad \leq \IP \Bigl( \sup_{\tau_1, \tau_2 \in[0,1]} \mathop{\sup_{\llvert  u
- \tau_1\rrvert   \leq\delta_n }}_{{\llvert  v - \tau_2\rrvert   \leq\delta_n}}
\bigl\llvert \hat G_{n, U} ( u, v; \omega) - \hat G_{n,U} (
\tau_1, \tau_2; \omega) \bigr\rrvert > x
\bigl((nb_n)^{-1/2}+b_n^k\bigr),
\nonumber
\\
&&\quad\hspace*{22pt}\sup_{\tau\in[0,1]} \bigl\llvert \hat F_n^{-1}(
\tau) - \tau\bigr\rrvert \leq \delta_n \Bigr) + \IP \Bigl(\sup
_{\tau\in[0,1]} \bigl\llvert \hat F_n^{-1}(\tau) -
\tau\bigr\rrvert > \delta_n \Bigr)
\nonumber
\\
&&\quad = P^n_1 + P^n_2,\qquad\mbox{say}.
\nonumber
\end{eqnarray*}
It follows from Lemma~\ref{lemquantproc} that $P^n_2$ is $\mathrm{o}(1)$. As
for $P^n_1$, it is bounded by
\begin{eqnarray}
&& \IP \Bigl( \sup_{\tau_1, \tau_2 \in[0,1]} \mathop{\sup_{\llvert  u -
\tau_1\rrvert   \leq\delta_n}}_{\llvert  v - \tau_2\rrvert   \leq\delta_n}
\bigl\llvert \hat H_{n,
U} ( u, v; \omega) - \hat H_{n,U} (
\tau_1, \tau_2; \omega) \bigr\rrvert >
\bigl(1+(nb_n)^{1/2}b_n^k\bigr)x/2
\Bigr)
\nonumber
\\
&&\quad{}+ I \Bigl\{ \sup_{\tau_1, \tau_2 \in[0,1]} \mathop{\sup
_{\llvert  u - \tau_1\rrvert   \leq\delta_n }}_{{ \llvert  v - \tau_2\rrvert   \leq\delta
_n}} \bigl\llvert \IE\hat G_{n, U}
( u, v; \omega) - \IE\hat G_{n,U} (\tau_1,
\tau_2; \omega) \bigr\rrvert > \bigl((nb_n)^{-1/2}+b_n^k
\bigr) x/2 \Bigr\},\nonumber
\end{eqnarray}
where
the first term tends to zero in view of~(\ref
{thmAsympDensityEstimatoreqnstochequicont}). To see that the
indicator in the second term also is $\mathrm{o}(1)$, note that
\begin{eqnarray*}
&& \sup_{\tau_1, \tau_2 \in[0,1]} \mathop{\sup_{\llvert  u - \tau_1\rrvert
\leq\delta_n }}_{\llvert  v - \tau_2\rrvert   \leq\delta_n}
\bigl\llvert \IE\hat G_{n, U} ( u, v; \omega) - \IE\hat
G_{n,U} (\tau_1, \tau_2; \omega) \bigr\rrvert
\\
&&\quad \leq \sup_{\tau_1, \tau_2 \in[0,1]} \mathop{\sup_{\llvert  u - \tau
_1\rrvert   \leq\delta_n }}_{ \llvert  v - \tau_2\rrvert   \leq\delta_n}
\bigl\llvert \IE\hat G_{n,
U} ( u, v; \omega) - \mathfrak{f}_{q_u, q_v}(
\omega) - B_n^{(k)}(u, v, \omega)\bigr\rrvert
\\
&&\qquad{} + \sup_{\tau_1, \tau_2 \in[0,1]} \mathop{\sup_{\llvert  u - \tau_1\rrvert
\leq\delta_n }}_{ \llvert  v - \tau_2\rrvert   \leq\delta_n}
\bigl\llvert B_n^{(k)}(\tau_1,
\tau_2, \omega) + \mathfrak{f}_{q_{\tau_1}, q_{\tau_2}}(\omega) - \IE\hat
G_{n,U} (\tau_1, \tau_2; \omega) \bigr\rrvert
\\
&&\qquad{} + \sup_{\tau_1, \tau_2 \in[0,1]} \mathop{\sup_{\llvert  u - \tau_1\rrvert
\leq\delta_n }}_{\llvert  v - \tau_2\rrvert   \leq\delta_n}
\bigl\llvert \mathfrak{f}_{q_u,
q_v}(\omega)+B_n^{(k)}(u,
v, \omega) - \mathfrak{f}_{q_{\tau_1},
q_{\tau_2}}(\omega) - B_n^{(k)}(
\tau_1, \tau_2, \omega) \bigr\rrvert
\\
&&\quad = \mathrm{o}\bigl( n^{-1/2} b_n^{-1/2} +
b_n^k\bigr) + \mathrm{O}\bigl(\delta_n\bigl(1+\llvert
\log\delta_n\rrvert \bigr)^{d}\bigr),
\end{eqnarray*}
where $d$ still is the constant from Lemma~\ref
{lemLipschitzLaplaceSD} corresponding to $j=k$. Here, we have applied
part (ii) of Theorem~\ref{thmAsympDensityEstimator} to bound the
first two terms and Lemma~\ref{lemLipschitzLaplaceSD} for the third
one. For any fixed $\omega$, thus, $P^n(\omega)=\mathrm{o}(1)$, which
establishes the pointwise version of the claim.

We now turn to the uniformity (with respect to $\omega$) issue. For an
arbitrary~$y_n > 0$, similar arguments as above yield, with the same
$\delta_n$,
\begin{eqnarray*}
&&\IP \Bigl( \sup_{\omega\in\IR}\sup_{\tau_1, \tau_2 \in[0,1]} \bigl
\llvert \hat G_{n,R} (\tau_1, \tau_2; \omega) -
\hat G_{n,U} (\tau_1, \tau _2; \omega) \bigr
\rrvert > y_n \Bigr)
\\
&&\quad \leq\IP \Bigl( \sup_{\omega\in\IR}\sup_{\tau_1, \tau_2 \in
[0,1]}
\mathop{\sup_{\llvert  u - \tau_1\rrvert   \leq\delta_n }}_{ \llvert  v - \tau_2\rrvert
\leq\delta_n} \bigl\llvert \hat
H_{n, U} ( u, v; \omega) - \hat H_{n,U} (\tau _1,
\tau_2; \omega) \bigr\rrvert > (nb_n)^{1/2}y_n/2
\Bigr)
\\
&&\qquad{} + I \Bigl\{\sup_{\omega\in\IR} \sup_{\tau_1, \tau_2 \in[0,1]}
\mathop{\sup_{\llvert  u - \tau_1\rrvert   \leq\delta_n }}_{ \llvert  v - \tau_2\rrvert   \leq
\delta_n} \bigl\llvert \IE\hat
G_{n, U} ( u, v; \omega) - \IE\hat G_{n,U} (
\tau_1, \tau_2; \omega) \bigr\rrvert > y_n/2
\Bigr\} + \mathrm{o}(1).
\end{eqnarray*}
That the indicator in the latter expression is $\mathrm{o}(1)$ follows by the
same arguments as above [note that Lemma~\ref{lemLipschitzLaplaceSD}
and the statement of part (ii) both hold uniformly in $\omega\in\IR
$]. To\vspace*{1pt} bound the probability term, observe that by Lemma~\ref
{lemBoundDFT}, $\sup_{\tau_1,\tau_2}\sup_{j=1,\ldots,n}
\llvert  I_{n,U}^{\tau_1,\tau_2}(2\uppi j/n)\rrvert   $ is~$ \mathrm{O}_P(n^{2/K})$ for any
$K>0$. Moreover, the uniform Lipschitz continuity of~$W$ implies that
$W_n$ also is uniformly Lipschitz continuous with constant of order
$\mathrm{O}(b_n^{-2})$. Combining those facts with Lemma~\ref
{lemLipschitzLaplaceSD} and the assumptions on~$b_n$, we obtain
\[
\mathop{\sup_{\omega_1,\omega_2 \in\IR}}_{{\llvert  \omega_1-\omega
_2\rrvert  \leq n^{-3}}} \sup_{\tau_1, \tau_2 \in[0,1]}
\bigl\llvert \hat H_{n, U} ( \tau_1,\tau_2;
\omega_1) - \hat H_{n,U} (\tau_1,
\tau_2; \omega _2) \bigr\rrvert = \mathrm{o}_P(1).
\]
By periodicity of $\hat H_{n, U}$ in the argument $\omega$, it thus
remains to show that
\[
\max_{\omega=0,2\uppi n^{-3},\ldots,2\uppi}\sup_{\tau_1, \tau_2 \in
[0,1]}\mathop{\sup
_{\llvert  u - \tau_1\rrvert   \leq\delta_n }}_{{ \llvert  v - \tau_2\rrvert
\leq\delta_n}} \bigl\llvert \hat H_{n, U} (
u, v; \omega) - \hat H_{n,U} (\tau _1, \tau_2;
\omega) \bigr\rrvert = \mathrm{o}_P(1).
\]
Lemmas \ref{lemThm224VW} and~\ref{lemOrderSmallIncrements} entail
the existence of a random variable $S(\omega)$ such that, for any
fixed $\omega\in\IR$,
\begin{eqnarray*}
\sup_{\tau_1, \tau_2 \in[0,1]}\mathop{\sup_{\llvert  u - \tau_1\rrvert   \leq
\delta_n }}_{{ \llvert  v - \tau_2\rrvert   \leq\delta_n}}
\bigl\llvert \hat H_{n, U} ( u, v; \omega) - \hat H_{n,U} (
\tau_1, \tau_2; \omega) \bigr\rrvert \leq\bigl\llvert S(
\omega)\bigr\rrvert + R_n(\omega),
\end{eqnarray*}
where $\sup_{\omega\in\IR} \llvert  R_n(\omega)\rrvert   = \mathrm{o}_P(1)$ and
\[
\max_{\omega=0,2\uppi n^{-3},\ldots,2\uppi}\IE\bigl[\bigl\llvert S^{2L}(\omega)
\bigr\rrvert \bigr] \leq K_L^{2L} \biggl(\int
_{0}^{\eta} \epsilon^{-4/(2L\gamma)} \,\mathrm{d}
\epsilon+ \bigl(\delta_n^{\gamma/2} + 2(nb_n)^{-1/2}
\bigr) \eta ^{-8/(2L\gamma)} \biggr)^{2L}
\]
for any $0< \gamma< 1, L \in\IN$, $0<\eta<\delta_n$, and a
constant $K_L$ depending on $L$ only.
For appropriate choice of $L$ and $\gamma$, this latter bound is
$\mathrm{o}(n^{-3})$; since the maximum is over a set with $\mathrm{O}(n^3)$ elements.
This completes the proof of part (iii).

\subsection{Details for the proof of parts \textup{(i)} and \textup{(iii)} of Theorem~\texorpdfstring{\protect\ref{thmAsympDensityEstimator}}{3.6}}\label{app-B}

This section contains the main lemmas used in Sections~\ref{Parti}
and~\ref{Partiii} above. We use the notation introduced at the
beginning of the proof of Theorem~\ref{thmAsympDensityEstimator}. The
proofs of the results presented here can be found in the
online supplement~\cite{supp} [Section~1.3].

For the statement of the first result, recall that, for any
non-decreasing, convex function $\Psi\dvt  \IR^+ \to\IR^+$ with $\Psi
(0) = 0$ the \textit{Orlicz norm} of a real-valued random variable $Z$
is defined as (see, e.g., van~der Vaart and Wellner \cite{vanderVaartWellner1996}, Chapter~2.2)
\[
\llVert Z\rrVert _\Psi= \inf\bigl\{C>0\dvt  \IE\Psi \bigl( {\llvert Z
\rrvert }/{C} \bigr) \leq 1\bigr\}.
\]

%
%
\begin{lemma}
\label{lemThm224VW}
Let $\{\IG_t\dvt  t \in T\}$ be a separable stochastic process with $\llVert
\IG_s - \IG_t\rrVert  _{\Psi} \leq C \,\dd (s,t)$ for all $s, t$ with $d(s,t)
\geq\bar\eta/2 \geq0$. Denote by $D(\epsilon,d)$ the packing
number of the metric space $(T,d)$. Then, for any $\delta>0 $, $\eta
\geq\bar\eta$, there exists a random variable $S_1$ and a constant
$K<\infty$ such that
\begin{eqnarray*}
\sup_{d(s,t) \leq\delta} \llvert \IG_s - \IG_t
\rrvert &\leq& S_1 + 2 \sup_{d(s,t) \leq\bar\eta, t \in\tilde T} \llvert
\IG_s - \IG _t\rrvert
\end{eqnarray*}
and
\begin{eqnarray*}
\llVert S_1\rrVert _\Psi&\leq& K \biggl[ \int
_{\bar\eta/2}^{\eta} \Psi^{-1} \bigl( D(\epsilon,
d) \bigr) \,\mathrm{d}\epsilon+ (\delta+ 2 \bar \eta ) \Psi^{-1} \bigl(
D^2(\eta, d) \bigr) \biggr],
\end{eqnarray*}
where the set $\tilde T$ contains at most $D(\bar\eta,d)$ points. In
particular, by Markov's inequality (cf. van~der Vaart and Wellner \cite{vanderVaartWellner1996}, page~96),
\begin{eqnarray*}
\IP \bigl( \llvert S_1\rrvert > x \bigr) \leq \biggl(\Psi \biggl( x
\biggl[ 8 K \biggl( \int_{\bar\eta/2}^{\eta}
\Psi^{-1} \bigl(D(\epsilon, d) \bigr) \,\mathrm{d}\epsilon+ (\delta+ 2
\bar\eta) \Psi^{-1} \bigl(D^2(\eta, d) \bigr) \biggr)
\biggr]^{-1} \biggr) \biggr)^{-1}
\end{eqnarray*}
for any $x>0$.
\end{lemma}

%
%
\begin{lemma} \label{lemSixthMomIncrementHn}
Let $X_0,\ldots,X_{n-1}$ be the finite realization of a strictly
stationary process with $X_0 \sim U[0,1]$, and let \textup{(W)} hold. For $x =
(x_1, x_2)$ let~$\hat H_n(x;\omega):= \sqrt{n b_n} (\hat
G_n(x_1,x_2;\omega) - \IE[\hat G_n(x_1,x_2;\omega)])$. For any Borel
set~$A$, define
\[
d^A_n(\omega):= \sum_{t=0}^{n-1}
I\{X_t \in A\} \ee^{-\ii t \omega}.
\]
Assume that, for $p=1,\ldots,P$, there exist a constant $C$ and a
function $g\dvtx  \IR^+ \to\IR^+$, both independent of $\omega
_1,\ldots,\omega_p \in\IR, n$ and $A_1,\ldots,A_p$, such that
\begin{equation}
\label{lemSixthMomIncrementHneqnAssumption} \bigl\llvert \cum\bigl(d_n^{A_1}(
\omega_1), \ldots, d_n^{A_p}(
\omega_p)\bigr) \bigr\rrvert \leq C \Biggl( \Biggl\llvert
\Delta_n \Biggl(\sum_{i=1}^p
\omega_i \Biggr) \Biggr\rrvert + 1 \Biggr) g(\eps)
\end{equation}
for any Borel sets $A_1, \ldots, A_p$ with $\min_j \IP(X_0 \in A_j)
\leq\varepsilon$. Then there exists a constant $K$ (depending on
$C,L,g$ only) such that
\[
\sup_{\omega\in\IR}\sup_{\llVert   a - b \rrVert  _1 \leq\varepsilon} \IE\bigl\llvert \hat
H_n(a;\omega) - \hat H_n(b;\omega) \bigr\rrvert
^{2L} \leq K \sum_{\ell=0}^{L-1}
\frac{g^{L-\ell}(\eps)}{(n b_n)^\ell}
\]
for all $\varepsilon$ with $g(\varepsilon)<1$ and all $L = 1, \ldots, P$.
\end{lemma}

%
%
\begin{lemma} \label{lemLipschitzLaplaceSD}
Under\vspace*{1pt} the assumptions of Theorem~\ref{thmAsympDensityRankEstimator},
the derivative $(\tau_1, \tau_2) \mapsto
\frac{\mathrm{d}^j}{\mathrm{d}\omega^j}\mathfrak{f}_{q_{\tau_1},
q_{\tau_2}}(\omega)$ exists and satisfies, for any $j \in\IN_0$
and some constants $C,d$ that are independent of
$a=(a_1,a_2),b=(b_1,b_2)$ but may depend on $j$,
\[
\sup_{\omega\in\IR} \biggl\llvert \frac{\mathrm{d}^j}{\mathrm{d}\omega
^j}
\mathfrak{f}_{q_{a_1},q_{a_2}}(\omega)-\frac{\mathrm{d}^j}{\mathrm{d}\omega^j}\mathfrak{f}_{q_{b_1}, q_{b_2}}(
\omega )\biggr\rrvert \leq C \llVert a-b\rrVert _1 \bigl(1+\bigl
\llvert \log\llVert a-b\rrVert _1\bigr\rrvert \bigr)^d.
\]
\end{lemma}

%
%
\begin{lemma} \label{lemOrderDftYinA}
Let the strictly stationary process $(X_t)_{t \in\IZ}$ satisfy assumption~\textup{(C)}.
For any Borel set $A$, define
\[
d^A_n(\omega):= \sum_{t=0}^{n-1}
I\{X_t \in A\} \ee^{-\ii t \omega}.
\]
Let $A_1, \ldots, A_p \subset[0,1]$ be intervals, and let $\eps:=
\min_{j=1,\ldots,p} \IP(X_0 \in A_j)$. Then, for any $p$-tuple
$\omega
_1,\ldots,\omega_p \in\IR$,
\[
\bigl\llvert \cum\bigl(d_n^{A_1}(\omega_1),
\ldots, d_n^{A_p}(\omega_p)\bigr) \bigr\rrvert
\leq C \Biggl( \Biggl\llvert \Delta_n \Biggl(\sum
_{i=1}^p \omega_i \Biggr) \Biggr
\rrvert + 1 \Biggr) \varepsilon\bigl(\llvert \log\eps\rrvert +1
\bigr)^d,
\]
where $\Delta_n(\lambda):= \sum_{t=0}^{n-1}\ee^{\ii t\lambda}$ and
the constants $C,d$ depend only on $K,p$, and~$\rho$ [with $\rho$
from condition~\textup{(C)}].
\end{lemma}

%
%
\begin{lemma}\label{lemquantproc}
Let $X_0,\ldots,X_{n-1}$ be the finite realization of a strictly
stationary process satisfying~\textup{(C)} and such that $X_0 \sim U[0,1]$. Then
\[
\sup_{\tau\in[0,1]} \bigl\llvert \hat F_n^{-1}(
\tau) - \tau\bigr\rrvert = \mathrm{O}_P\bigl(n^{-1/2}\bigr).
\]
\end{lemma}

%
%
\begin{lemma}\label{lemBoundDFT}
Let the strictly stationary process $(X_t)_{t \in\IZ}$ satisfy
assumption~\textup{(C)}; assume moreover that $X_0 \sim U[0,1]$. For any $y\in
[0,1]$, define
\[
d_n^{y}(\omega):= \sum_{t=0}^{n-1}
I\{X_t \leq y\}\ee^{-\ii\omega t}.
\]
Then, for any $k\in\IN$,
\[
\sup_{\omega\in\cF_n} \sup_{y \in[0,1]} \bigl\llvert
d_n^{y}(\omega) \bigr\rrvert = \mathrm{O}_P
\bigl(n^{1/2+1/k}\bigr).
\]
\end{lemma}

%
%
\begin{lemma}\label{lemOrderSmallIncrements}
Under the assumptions of Theorem~\ref{thmAsympDensityEstimator}, let
$\delta_n$ be a sequence of non-negative real numbers. Assume that
there exists $\gamma\in(0,1)$, such that~$\delta_n =\mathrm{O}( (n
b_n)^{-1/\gamma})$. Then
\[
\sup_{\omega\in\IR}\mathop{\sup_{u,v \in[0,1]^2 }}_{{ \llVert   u - v \rrVert
_1 \leq\delta_n }}
\bigl\llvert \hat H_n(u;\omega) - \hat H_n(v;\omega)
\bigr\rrvert = \mathrm{o}_P(1).
\]
\end{lemma}
\end{appendix}

\section*{Acknowledgements}
The authors are grateful to three referees
for their constructive comments on an earlier version of this paper, which
led to a substantial improvement of our work.

This work has been supported by the Sonderforschungsbereich
``Statistical modelling of
nonlinear dynamic processes'' (SFB~823) of the Deutsche Forschungsgemeinschaft.
Tobias Kley was supported by a PhD Grant of the Ruhr-Universit\"at
Bochum and by the
Ruhr-Universit\" at Research School funded by Germany's Excellence
Initiative [DFG GSC 98/1].
Stanislav Volgushev and Holger Dette were supported by the
Collaborative Research Center ``Statistical
modeling of nonlinear dynamic processes'' (SFB 823, Teilprojekt
A1, C1) of the German Research Foundation (DFG).
Marc Hallin was supported by the Belgian Science Policy Office (2012--2017)
Interuniversity Attraction Poles and a Humboldt-Forschungspreis of the
Alexander von Humboldt-Stiftung.

\begin{supplement}
\stitle{Supplement to ``Quantile spectral processes: Asymptotic~analysis and inference''.}
\slink[doi]{10.3150/15-BEJ711SUPP} 
\sdatatype{.pdf}
\sfilename{BEJ711\_supp.pdf}
\sdescription{We provide details for the proof of part (ii) of
Theorem \ref{thmAsympDensityEstimator}, and proofs for Propositions
\ref{propalph}, \ref{propgmc}, and \ref{propinr}. Further, we prove results from Section \ref{app-B},
namely Lemmas \ref{lemThm224VW}--\ref{lemOrderSmallIncrements}.}
\end{supplement}

%

\printhistory
\end{document}